\documentclass[trsc,nonblindrev]{informs3} 

\OneAndAHalfSpacedXI 

\usepackage{afterpage}
\usepackage{caption}
\usepackage{amssymb}
\usepackage{float}
\usepackage{dcolumn}
\usepackage{multirow}
\usepackage{color}
\usepackage{enumitem}
\usepackage{slashbox}
\usepackage{pict2e}
\usepackage{microtype}
\usepackage[noend]{algorithmic}
\usepackage{algorithm}
\usepackage{rotating}
\usepackage{pdflscape}
\usepackage[table]{xcolor}
\usepackage{tabularx, booktabs}
\usepackage{algorithmic}
\usepackage{rotating} 
\newcolumntype{M}[1]{>{\centering}m{#1}}
\newcolumntype{F}{>{\centering\arraybackslash}m{1.5cm}}

\usepackage{natbib}
 \bibpunct[, ]{(}{)}{,}{a}{}{,}%
 \def\bibfont{\small}%
\TheoremsNumberedThrough     

\EquationsNumberedThrough    

\begin{document}


\RUNAUTHOR{Matl, Hartl, and Vidal}

\RUNTITLE{Workload Equity in Vehicle Routing Problems: A Survey and Analysis}

\TITLE{Workload Equity in Vehicle Routing Problems: \protect \\ A Survey and Analysis}

\ARTICLEAUTHORS{%
\AUTHOR{P.~Matl, R.F.~Hartl}
\AFF{University of Vienna, Austria, \EMAIL{piotr.matl@univie.ac.at}, \EMAIL{richard.hartl@univie.ac.at}\URL{}}
\AUTHOR{T.~Vidal}
\AFF{Pontifical Catholic University of Rio de Janeiro, Brazil, \EMAIL{vidalt@inf.puc-rio.br}, \URL{}}
} 

\ABSTRACT{%
Over the past two decades, equity aspects have been considered in a growing number of models and methods for vehicle routing problems (VRPs). Equity concerns most often relate to fairly allocating workloads and to balancing the utilization of resources, and many practical applications have been reported in the literature. However, there has been only limited discussion about how workload equity should be modelled in the context of VRPs, and various measures for optimizing such objectives have been proposed and implemented without a critical evaluation of their respective merits and consequences. 
\newline
\indent This article addresses this gap by providing an analysis of classical and alternative equity functions for bi-objective VRP models. In our survey, we review and categorize the existing literature on equitable VRPs. In the analysis, we identify a set of axiomatic properties which an ideal equity measure should satisfy, collect six common measures of equity, and point out important connections between their properties and the properties of the resulting Pareto-optimal solutions. To gauge the extent of these implications, we also conduct a numerical study on small bi-objective VRP instances solvable to optimality.
Our study reveals two undesirable consequences when optimizing equity with non-monotonic functions: Pareto-optimal solutions can consist of non-TSP-optimal tours, and even if all tours are TSP-optimal, Pareto-optimal solutions can be workload inconsistent, i.e.~composed of tours whose workloads are all equal to or longer than those of other Pareto-optimal solutions. We show that the extent of these phenomena should not be under-estimated. The results of our bi-objective analysis remain valid also for weighted sum, constraint-based, or single-objective models. Based on this analysis, we conclude that monotonic equity functions are more appropriate for certain types of VRP models, and suggest promising avenues for further research on equity in logistics.
}%



\maketitle

%


\section{Introduction.}
\label{1-intro}

Logistics is multi-objective -- the complexity of real-life logistics planning and decision-making often cannot be reduced to cost only. From a more holistic point of view, non-monetary factors such as service quality, reliability, customer satisfaction, and consistency can be decisive for an organization's medium or long-term performance, especially if those factors can be improved at little additional cost (\citealt{jozefowiez2008}, \citealt{kovacs2014}).

As a consequence, concerns about workload equity have gained increasing attention during the past two decades. This is not surprising, since the explicit consideration of equity issues can provide tangible benefits to organizations. For example, balancing the capacity utilization of vehicles can free up bottleneck resources for future growth in operations without the need for additional investment \citep{apte2006}. Equitable workload allocation can also improve the acceptance of operational plans among drivers, their morale, and the quality of customer service provided \citep{liu2006}. The relevance of these concerns is reflected in the number of practical applications found in the operations research and logistics literature \citep{karsu2015}.

Despite this broad interest, there has been only limited discussion about how to model workload equity in the context of vehicle routing problems. Various equity objectives have been proposed and implemented, but the consequences of any specific equity objective choice remain largely unclear. It appears that this important modelling choice is often decided on an ad hoc basis. As a result, there is a lack of general guidelines for which \textit{types} of measures are appropriate in which contexts. Yet such guidance is essential given the broad variety of surveyed applications. Furthermore, paradoxes such as ``artificial equity'' have been observed empirically \citep{jozefowiez2002, halvorsen2016}, but remain unexplained in the general case.

Although single-objective comparative analyses of some equity objectives have been conducted \citep{campbell2008, huang2012, bertazzi2015, halvorsen2016}, the conclusions of these studies are of limited scope in more general bi-objective settings. Indeed, it is not obvious to what extent single-objective worst-case comparisons reflect the properties of intermediate trade-off solutions. Yet in practical applications, the primary interest lies ultimately in compromise solutions rather than any single-objective optimum.

In this survey and analysis, we take a step back and re-evaluate different ways of accommodating workload equity in bi-objective VRP models. Our contributions are fourfold:

\vspace*{0.3cm}
\begin{itemize}
	\item we review and categorize the literature on models and methods for equitable VRPs (Section~\ref{2-review}),
	\item we analyse the choice of equity measure from an axiomatic perspective and compare six common measures (Sections~\ref{assumptions} to \ref{measures}),
	\item we outline relevant connections between the properties of the chosen equity measure and the properties of Pareto-optimal solution sets (Section~\ref{implications}), and
	\item we conduct a numerical study to estimate the extent of some of these connections and to provide further insight into other relevant aspects for choosing an equity measure (Section~\ref{4-numerical}).
\end{itemize}
\vspace*{0.3cm}

We conclude this article in Section~\ref{5-conclusion} by summarizing our findings and suggesting a number of relevant directions for research on equity aspects in logistics applications.


\newpage
\section{Models and Methods for Equity in Vehicle Routing.}
\label{2-review}

Table~\ref{surveyTable} classifies the current literature on vehicle routing with equity considerations. Based on our search, this literature can be divided into five groups. A first stream of research focuses on a prototypical bi-objective capacitated VRP (CVRP) model with equity as the second objective (called the vehicle routing problem with route balancing -- VRPRB), and a second growing group of papers proposes time window extensions of this basic model. A third set of publications deal with the single-objective ``min-max VRP'', in which equity is the primary optimization objective. Some other works focus on incorporating equity issues into more complex VRP variants such as periodic VRPs or arc routing problems. These publications build the fourth stream. Finally, a large and varied selection of VRP application papers forms the fifth group.


\begin{table}[htbp]
\caption{Vehicle Routing Publications Incorporating Equity Concerns}
\centering
\scalebox{0.75}
{
\begin{tabular}{@{}lccccccccccc@{}}
\toprule
\multirow{2}{*}{{\bf Publication (by year)}} & \multicolumn{3}{c}{{\bf Equity Function}} & \multicolumn{2}{c}{{\bf Equity Metric}} & \multicolumn{4}{c}{{\bf Opt.~Model}} & \multicolumn{2}{c}{{\bf Opt.~Method}}  \\ \cmidrule(l){2-12} 
                                             & range     & min-max    & other    & length    & demand    & MO       & WS       & CN         & PO       & Exact    & Heur.  \\ \midrule
{\bf Standard VRPRB}               					 &           &            &          &           &           &          &          &            &          &          & $\times$ \\ \midrule
~~~\cite{lacomme2015}										     & $\times$  &            &          & $\times$  &           & $\times$ &          &            &          &          & $\times$ \\
~~~\cite{oyola2014}                          & $\times$  &            &          & $\times$  &           & $\times$ &          &            &          & $\times$ &        \\
~~~\cite{sarpong2013}                        &           & $\times$   &          & $\times$  & $\times$  & $\times$ &          &            &          & $\times$ &        \\
~~~\cite{reiter2012}                         &           & $\times$   &          & $\times$  &           & $\times$ &          &            &          &          & $\times$ \\
~~~\cite{jozefowiez2009}                     & $\times$  &            &          & $\times$  &           & $\times$ &          &            &          &          & $\times$ \\
~~~\cite{borgulya2008}                       & $\times$  &            &          & $\times$  &           & $\times$ &          &            &          &          & $\times$ \\
~~~\cite{jozefowiez2007}                     & $\times$  &            &          & $\times$  &           & $\times$ &          &            &          &          & $\times$ \\
~~~\cite{pasia2007a}                         & $\times$  &            &          & $\times$  &           & $\times$ &          &            &          &          & $\times$ \\
~~~\cite{pasia2007b}                         & $\times$  &            &          & $\times$  &           & $\times$ &          &            &          &          & $\times$ \\
~~~\cite{jozefowiez2006}                     & $\times$  &            &          & $\times$  &           & $\times$ &          &            &          &          & $\times$ \\
~~~\cite{jozefowiez2002}                     & $\times$  &            &          & $\times$  &           & $\times$ &          &            &          &          & $\times$ \\ \midrule
{\bf VRPRB with Time Windows}                &           &            &          &           &           &          &          &            &          &          &        \\ \midrule
~~~\cite{deFreitas2014}                      & $\times$  &            &          & $\times$  &           & $\times$ &          &            &          &          & $\times$ \\
~~~\cite{melian2014}                         & $\times$  &            &          & $\times$  &           & $\times$ &          &            &          &          & $\times$ \\
~~~\cite{banos2013}                          & $\times$  &            &          & $\times$  & $\times$  & $\times$ &          &            &          &          & $\times$ \\ \midrule
{\bf Min-Max VRP}                            &           &            &          &           &           &          &          &            &          &          &        \\ \midrule
~~~\cite{bertazzi2015}                       & $\times$  & $\times$   &          & $\times$  &           &          &          &            & $\times$ & $\times$ &        \\
~~~\cite{wang2014}                           &           & $\times$   &          & $\times$  &           &          &          &            & $\times$ &          & $\times$ \\
~~~\cite{narasimha2013}                      &           & $\times$   &          & $\times$  &           &          &          &            & $\times$ &          & $\times$ \\
~~~\cite{schwarze2013}                       &           & $\times$   &          & $\times$  &           &          &          &            & $\times$ &          & $\times$ \\
~~~\cite{yakici2013}                         &           & $\times$   &          & $\times$  &           &          &          &            & $\times$ &          & $\times$ \\
~~~\cite{huang2012}                          & $\times$  & $\times$   & $\times$ & $\times$  & $\times$  &          & $\times$ &            &          & $\times$ &        \\
~~~\cite{valle2011}                          &           & $\times$   &          & $\times$  &           &          &          &            & $\times$ & $\times$ & $\times$ \\
~~~\cite{carlsson2009}                       &           & $\times$   &          & $\times$  &           &          &          &            & $\times$ &          & $\times$ \\
~~~\cite{campbell2008}                       & $\times$  & $\times$   &          & $\times$  &           &          &          &            & $\times$ & $\times$ & $\times$       \\
~~~\cite{saliba2006}                         &           & $\times$   &          & $\times$  &           &          &          &            & $\times$ &          & $\times$ \\
~~~\cite{applegate2002}                      &           & $\times$   &          & $\times$  &           &          &          &            & $\times$ & $\times$ &        \\
~~~\cite{golden1997}                         &           & $\times$   &          & $\times$  &           &          &          &            & $\times$ &          & $\times$ \\ \midrule
{\bf Other VRPs with Equity Aspects}         &           &            &          &           &           &          &          &            &          &          &        \\ \midrule
~~~\cite{halvorsen2016}                      & $\times$  & $\times$   & $\times$ & $\times$  &           & $\times$ &          &            &          & $\times$ &        \\
~~~\cite{mandal2015}                         & $\times$  &            &          & $\times$  &           & $\times$ &          &            &          &          & $\times$ \\
~~~\cite{martinez2014}                       & $\times$  &            &          & $\times$  &           & $\times$ &          &            &          &          & $\times$ \\
~~~\cite{liu2013}                            &           & $\times$   &          &  $\times$ &           &          &          &            & $\times$ &          & $\times$ \\
~~~\cite{gulczynski2011}                     & $\times$  &            &          &           & $\times$  &          & $\times$ &            &          &          & $\times$ \\
~~~\cite{mei2011}                            &           & $\times$   &          & $\times$  &           & $\times$ &          &            &          &          & $\times$ \\
~~~\cite{mourgaya2007}                       &           & $\times$   &          &           & $\times$  &          &          & $\times$   &          & $\times$ &        \\
~~~\cite{lacomme2006}                        &           & $\times$   &          & $\times$  &           & $\times$ &          &            &          &          & $\times$ \\
~~~\cite{ribeiro2001}                        &           &            & $\times$ &           & $\times$  &          & $\times$ &            &          & $\times$ &        \\ \midrule
{\bf Application Papers}                     &           &            &          &           &           &          &          &            &          &          &        \\ \midrule
~~~\cite{deArmas2015}                        & $\times$  &            &          & $\times$  &           &          &          & $\times$   &          &          & $\times$ \\
~~~\cite{goodson2014}                        & $\times$  &            &          & $\times$  &           &          &          & $\times$   &          &          & $\times$ \\
~~~\cite{lopez2014}                          &           & $\times$   &          & $\times$  &           &          &          &            & $\times$ &          & $\times$ \\
~~~\cite{bektas2013}                         & $\times$  &            &          & $\times$  &           &          &          & $\times$   &          & $\times$ &        \\
~~~\cite{perugia2011}                        &           &            & $\times$ & $\times$  &           &          &          & $\times$   &          &          & $\times$ \\
~~~\cite{rienthong2011}                      & $\times$  &            &          & $\times$  &           &          &          &            & $\times$ & $\times$ &        \\
~~~\cite{kritikos2010}                       &           &            & $\times$ &           & $\times$  &          & $\times$ &            &          &          & $\times$ \\
~~~\cite{groer2009}                          &           &            & $\times$ & $\times$  & $\times$  &          &          & $\times$   &          &          & $\times$ \\
~~~\cite{mendoza2009}                        &           &            & $\times$ &           & $\times$  &          &          & $\times$   &          &          & $\times$ \\
~~~\cite{parragh2009}                        &           &            & $\times$ & $\times$  &           & $\times$ &          &            &          &          & $\times$ \\
~~~\cite{melachrinoudis2007}                 &           &            & $\times$ & $\times$  &           &          & $\times$ &            &          &          & $\times$ \\
~~~\cite{apte2006}                           & $\times$  &            &          & $\times$  & $\times$  &          &          & $\times$   &          &          & $\times$ \\
~~~\cite{jang2006}                           &           &            & $\times$ & $\times$  &           &          & $\times$ &            &          &          & $\times$ \\
~~~\cite{kim2006}                            & $\times$  &            &          & $\times$  &           &          &          & $\times$   &          &          & $\times$ \\
~~~\cite{lin2006}                            & $\times$  &            &          & $\times$  & $\times$  & $\times$ &          &            &          &          & $\times$ \\
~~~\cite{liu2006}                            & $\times$  &            &          & $\times$  & $\times$  &          &          & $\times$   &          &          & $\times$ \\
~~~\cite{pacheco2006}                        &           & $\times$   &          & $\times$  &           & $\times$ &          &            &          &          & $\times$ \\
~~~\cite{spada2005}                          &           & $\times$   &          & $\times$  &           &          &          &            & $\times$ &          & $\times$ \\
~~~\cite{blakeley2003}                       &           &            & $\times$ &           & $\times$  &          & $\times$ & $\times$   &          &          & $\times$ \\
~~~\cite{corberan2002}                       &           & $\times$   &          & $\times$  &           & $\times$ &          &            &          &          & $\times$ \\
~~~\cite{li2002}                             &           &            & $\times$ & $\times$  & $\times$  &          &          & $\times$   &          &          & $\times$ \\
~~~\cite{colorni2001}                        &           & $\times$   &          & $\times$  &           &          &          &            & $\times$ &          & $\times$ \\
~~~\cite{serna2001}                          &           & $\times$   &          & $\times$  &           &          &          &            & $\times$ &          & $\times$ \\
~~~\cite{lee1999}                            &           &            & $\times$ & $\times$  &           &          & $\times$ &            &          &          & $\times$ \\
~~~\cite{bowerman1995}                       &           &            & $\times$ &           & $\times$  &          & $\times$ &            &          &          & $\times$ \\ \bottomrule
\midrule
\end{tabular}
}
\label{surveyTable}
\end{table}

In this article, we make a distinction between equity \textit{metrics}, equity \textit{functions}, and equity \textit{objectives}. With equity metrics we refer to the resource that is balanced, e.g.~workload can be measured by the tour length, the demand served, etc. An equity function specifies how an index value is calculated for a given allocation of workloads, e.g.~using the range or standard deviation. We use the terms equity function and equity measure interchangeably. Finally, an equity objective refers to a specific combination of metric and function, ~e.g. the equity objective in the VRPRB is the range of tour lengths.

Most works consider tour length as the equity metric, based either on distance or duration. Fewer papers consider workload in terms of the demand served per tour. The two most common functions used to measure equity at the solution level are the largest workload (min-max), and the difference between the largest and the smallest (the range). More sophisticated functions, such as the standard deviation or the average deviation to the mean, are used in a few applications.

Various optimization models are used to accommodate equity concerns. In single-objective models, equity is either the primary objective (PO), or modelled as a constraint~(CN). In bi-objective models, equity is explicitly an objective and the set of efficient solutions is generated either with a weighted sum approach (WS) or with multi-objective methods (MO). Many surveyed publications focus on heuristics, but exact algorithms and discussions of analytical properties have also been presented.

\subsection{The VRP with Route Balancing.}
\label{2-1}

Among the surveyed papers, one standard bi-objective model has gained traction in the ``theoretical literature'', namely the VRPRB. This model extends the basic capacitated VRP (CVRP) with a second objective for equity, measured by the range of tour lengths. A formal mathematical formulation is given in \cite{oyola2014}.

The VRPRB was first introduced in \cite{jozefowiez2002}, with a parallel genetic algorithm combined with tabu search. Further work by the same authors was published in \cite{jozefowiez2006} with adaptations to NSGA-II based on parallel search and enhanced management of population diversity. This was followed by a paper introducing target-aiming Pareto search, a hybridization of NSGA-II, tabu search, and goal programming techniques \citep{jozefowiez2007}.

A number of competing solution procedures have since been proposed by other authors. \cite{pasia2007a} combine a randomized Clarke-Wright savings algorithm with Pareto local search. This work is further extended in \cite{pasia2007b} with Pareto ant colony optimization. \cite{borgulya2008} proposes an evolutionary algorithm in which recombination is omitted in favour of an adaptive mutation operator. Finally, \cite{jozefowiez2009} present a classical evolutionary algorithm extended with parallel search and elitist diversification management.

In contrast to earlier approaches, more recent contributions do not apply evolutionary algorithms. \cite{oyola2014} present a GRASP heuristic combined with local search. 
They reimplement a sequential version of the algorithm proposed in \cite{jozefowiez2009}, yielding better solution quality, but at the expense of higher computation times. The most recent article on the VRPRB by \cite{lacomme2015} introduces a multi-start decoder-based method which alternates between direct and indirect solution spaces. Non-dominated sets of solutions are extracted from TSP giant tour representations through a bi-objective split procedure and used as starting points for further search by path re-linking.

To the best of our knowledge, \cite{reiter2012} present the only exact algorithm designed specifically for a bi-objective VRP with an equity objective. However, in this work equity is measured with the maximum tour length, rather than the range. The authors propose an adaptive \mbox{$\epsilon$-constraint} method based on a hybridization of branch-and-cut and two genetic algorithms. An interesting conclusion of their study is that computing higher equity solutions requires considerably more computational effort than for lower equity (and thus lower cost) solutions. A more general model is investigated by \cite{sarpong2013}. The authors consider bi-objective VRPs with general min-max objectives, and a column generation algorithm is proposed to solve these problems to optimality.

\subsection{Time Window Extensions to the VRPRB.}
\label{2-2}

Several publications consider the addition of time windows to the VRPRB. \cite{banos2013} examine the range equity objective based on either tour lengths or tour loads. 
Their results suggest that the length and load based equity objectives lead to noticeably different solutions.

A hybrid heuristic combining iterated local search, randomized variable neighbourhood descent, and genetic recombination is proposed in \cite{deFreitas2014}. The authors compare their approach with NSGA-II and a reimplementation of the algorithm in \cite{banos2013}, finding that their approach outperforms the others. 

\cite{melian2014} present a bi-objective scatter search heuristic for the VRPRB with time windows. The method is shown to approximate well the exact Pareto front on small instances, and outperforms NSGA-II on larger ones. The authors note that the addition of time windows leads to undesired consequences when equity is optimized: unnecessary waiting times may be added to tours in order to artificially improve the equity objective.

\subsection{The Min-Max VRP.}
\label{2-3}

In contrast to the provider-oriented perspective of minimizing total cost, service and client-oriented logistics applications often consider equity issues as the primary optimization objective. Examples include school bus routing (minimizing the longest time spent in transit in \citealt{serna2001}), delivery of newspapers or perishable goods (minimizing the latest delivery time in \citealt{applegate2002}), communication or computer networks (minimizing the maximum latency in \citealt{nace2008}), humanitarian relief (minimizing the latest arrival of aid in \citealt{golden2014}), as well as military reconnaissance or surveillance operations (minimizing the longest mission time in \citealt{carlsson2013}). A common feature of many of these applications is that the equity objective is defined with the min-max measure, and post-optimization procedures are applied to optimize the remaining parts of the solution if necessary.

An early contribution to the class of min-max VRP problems can be found in \cite{golden1997}. The authors apply a min-max length objective to the CVRP, the CVRP with multiple use of vehicles, and the $m$-TSP. Their computational experiments show that the total cost of the best min-max solution can be significantly higher than the cost of the best min-sum solution. \cite{applegate2002} revisit a min-max VRP instance from a competition based on newspaper delivery and prove the optimality of the heuristic solution which had won the original competition. \cite{saliba2006} presents construction heuristics for a lexicographic version of the min-max VRP.

Extensions to the min-max VRP have recently gained attention. The multi-depot variant (min-max MDVRP) was first introduced in \cite{carlsson2009}. Their paper presents a theoretical analysis of optimal solutions and compares three different heuristics. This work is further pursued in \cite{narasimha2013}, in which an ant colony optimization method is combined with equitable region partitioning to solve a series of single-depot problems. \cite{wang2014} also extend some of the ideas presented in \cite{carlsson2009} to the min-max MDVRP and compare the performance of four different heuristics.

More complex extensions of the single-depot case are presented in \cite{yakici2013} and \cite{schwarze2013}. \cite{yakici2013} consider a min-max VRP with multiple regions, multiple types of customer demand, and a heterogeneous fleet of vehicles, with the possibility of split deliveries. \cite{schwarze2013} study the Skill VRP, which arises in the field of service technician routing. Based on the observation that solutions to the standard model produce TSP-like solutions with poor capacity utilization of the technicians' skill levels, an alternative model with a min-max objective is introduced to better balance the utilization of these resources. Finally, \cite{valle2011} describe the min-max selective VRP, a variant encountered when designing wireless sensor networks. 

Since all of the above publications propose single-objective models, an important question is to what extent the chosen objective affects the properties and structure of optimal solutions. This topic is studied in \cite{campbell2008, huang2012}, and \cite{bertazzi2015}. These articles will be reviewed in Section~\ref{3-analysis}.


\subsection{Other VRPs with Equity Aspects.}
\label{2-4}

Some studies examine equity objectives from a multi-period perspective. \cite{gulczynski2011} consider two variants of the periodic VRP (PVRP), including one with an equity objective. An interesting aspect of their model is that the workload is measured by the number of customers per tour, which can be the main determinant of workload in certain applications. \cite{gulczynski2011} measure equity with the range of workloads on a given day and combine this with the classical total cost objective into a weighted sum objective function. \cite{mourgaya2007} examine the PVRP from a tactical perspective to optimize regional compactness and workload balance among routes. In their model, a route's workload is represented by the total demand served, and handled with a constraint on the maximum allowable workload.

In contrast, workload equity is the primary objective in the PVRP considered in \cite{liu2013}. The problem is motivated by home healthcare logistics, in which reasonable working hours are a major concern. The authors therefore minimize the maximum route duration as the primary objective. \cite{ribeiro2001} examine a PVRP with three objectives: minimization of total cost, minimization of the standard deviation of route durations, and minimization of an inconsistency objective for customer service. However, only a small instance is solved due to the non-linearity of the proposed model.

\cite{martinez2014} present a transportation location routing problem motivated by soft drink distribution. The strategic decisions of locating distribution centres are jointly optimized with the operational routing decisions. The authors formulate a bi-objective model that minimizes total system cost as well as the range of delivery route lengths. 

Several publications deal with bi-objective capacitated arc routing problems (CARPs) with equity objectives. \cite{lacomme2006} consider the trade-off between total distance minimization and minimization of the maximum tour cost, and propose an adaptation of NSGA-II for determining the set of efficient solutions. The bi-objective CARP is examined also in \cite{mei2011}, who propose a decomposition-based memetic algorithm. An extensive comparative study with the algorithm of \cite{lacomme2006} and a standard NSGA-II shows that the approach described in \cite{mei2011} yields significant improvements in solution quality, but at the expense of longer computation times.

A generalization of the CARP and CVRP -- the mixed capacitated general routing problem (MCGRP) -- is studied in \cite{mandal2015}. The authors measure equity with the range of route costs, and propose an extension of NSGA-II with Pareto-based local search. Computational experiments suggest that some crossover operators are more effective than others for optimizing the balance criterion. In \cite{halvorsen2016} small instances of the MCGRP are solved exactly with four different functions for the equity objective. We review the observations and conclusions of this study in more detail in Section~\ref{3-analysis}.

\subsection{Applications.}
\label{2-5}

Some of the earliest VRP applications with equity aspects focus on school bus routing. Since these transportation services are generally provided by the public sector, equity objectives must be considered in addition to cost efficiency. \cite{bowerman1995} consider the design of school bus routes in an urban setting in Canada. Next to minimizing the number of bus routes, equity objectives are included to balance the route lengths and the number of students served by the routes, as well as to minimize student walking distance. \cite{serna2001} re-examine a rural school bus routing problem in Spain from a social perspective. A min-max VRP model is optimized with the aim of minimizing the longest time spent in transit. Computational tests show that the min-max solution can significantly improve the service objective compared to the solution of the standard min-sum model, while still reducing costs compared to manual plans.


\cite{corberan2002} extend the single-objective min-max VRP of \cite{serna2001} to two objectives, and analyse explicitly the trade-off between minimizing the number of bus routes and minimizing the maximum travel time. This bi-objective model is further investigated in \cite{pacheco2006}, who propose new construction and tabu search heuristics to obtain improved solutions compared to those in \cite{corberan2002}.

Other school bus routing case studies include those of \cite{li2002} and \cite{spada2005}. \cite{li2002} describe the planning of bus services for a kindergarten in Hong Kong. They optimize hierarchically a set of four objectives, including one to balance the load and the travel time among the buses. \cite{spada2005} consider school bus routing with a focus on quality of service and present case studies of three towns in Switzerland. Service quality is operationalized as the child's time loss, measured as the difference between transit plus waiting time and the shortest possible travel time. \cite{spada2005} analyse the differences between optimizing the sum of time losses and optimizing the largest time loss.

Private sector bus services have also been studied with equity issues in mind. \cite{perugia2011} describe the design of a home-to-work bus service for employees of a large research center in Rome, Italy. A generalized bi-objective VRPTW is proposed to optimize the total cost as well as the quality of the service, measured as the passenger's time loss compared to direct travel. An interesting feature of this model is that equity in the distribution of time losses is enforced indirectly by appropriately setting time windows on the arrival times of buses at different stops. A similar case is studied by \cite{lopez2014}, who consider the planning of bus services for employees of a large company in Spain. The problem is modelled as an open VRP with a min-max objective in order to balance the time spent in transit and comply with insurance constraints. 

Dial-a-ride services, such as the door-to-door transportation of handicapped or elderly persons, are also a relevant field of application for equitable VRPs. Successful dial-a-ride operations must account for quality of service. For this reason, it is common to impose constraints on the time loss and/or the deviation to planned time windows to ensure that no user's experience is excessively worse compared to others' \citep{cordeau2007}. However, in some studies this aspect is also explicitly considered in the optimization objective(s) in order to actively favour service equity. For instance, \cite{colorni2001} investigate dial-a-ride services for two cities in Italy and examine the impact of minimizing the maximum time loss among all passengers. \cite{melachrinoudis2007} describe the case of client transportation for a health care organization in Boston, Massachusetts. They consider the trade-off between total cost and total client inconvenience, the latter based on excess ride time, earliness, and lateness. By generating Pareto sets on sample cases, the authors guide the decision-maker in selecting suitable weights for a single weighted sum objective used in their proposed decision support tool. \cite{parragh2009} consider recurring patient transportation services at the Austrian Red Cross. In addition to the standard constraint on maximum transit time, they formulate a bi-objective model to minimize the average transit time and the total cost. 

Equity concerns appear as objectives also in a variety of other applications. \cite{lin2006} examine a location routing problem faced by a telecommunications company in Hong Kong. The authors solve a bi-objective formulation with the aim of minimizing total cost as well as range equity objectives for working hours and for demand served. They also consider how allowing multiple trips per employee affects the optimization of the equity objectives. \cite{deArmas2015} present a rich VRP application for a logistics services provider in Spain. In addition to the classical objectives of minimizing total cost and the number of vehicles, the authors also consider range objectives for route duration and route distance. 

In other applications, equity concerns are handled with constraints. \cite{liu2006} consider the planning of third party logistics services to convenience stores in Taiwan. To reduce overtime and improve acceptance of the routing plan among drivers, the total travel distance is minimized subject to constraints on the duration and the load of the vehicles. \cite{goodson2014} describes the election day routing of rapid response attorneys to poll observers in the United States. The fair allocation of workload is particularly important in this context because the attorneys are volunteers. Two alternative models are examined: one with equity based on the number of assigned polling stations per attorney, and the other based on their route durations.

\cite{groer2009} present a balanced billing cycle VRP encountered by an utility company. When accounts are cancelled or opened, fixed meter reading routes must be adjusted periodically to maintain efficiency while meeting regulatory requirements and customer service considerations. At the same time, it is desirable to balance the workload per day to avoid fixed periodic overtime costs and peaks in administrative work. This is handled by lower and upper bound constraints on the meter route distance and the number of meters serviced per route. A similar application is described in \cite{blakeley2003} for a complex technician routing and scheduling problem. The authors optimize periodic maintenance operations of elevators at various customer sites by minimizing a weighted sum of total travel cost, overtime, and workload balance, the latter being based on customer service times. In addition, upper and lower bounds on the workload of each route are included, as it is noted that some weight combinations would otherwise lead to unacceptable day-to-day variations of workload.

Equitable VRPs arise also when optimizing delivery operations at libraries. \cite{apte2006} investigate the delivery of requested, returned, and new items between the branches of a large urban library network in San Francisco. Since the capacity of a delivery system depends critically on bottleneck resources, balancing the utilization of available truck capacities is essential for handling larger delivery volumes and leaving room for future growth in operations without the need for additional investments. In order to rebalance the delivery operations of the library network, the authors enforce lower and upper bounds on the number of stops and the load of each vehicle.

A similar application case is reported in \cite{rienthong2011} for a mobile library service on the Isle of Wight. The problem aims to determine a set of $m$ TSP routes for the mobile library such that every location is visited once during a recurring planning horizon of $m$ days. In order to balance the daily workload of the library, they add upper and lower bound constraints on the sum of each tour's travel and service time. An integer programming formulation of this problem is further investigated in \cite{bektas2013}.

Equity considerations can also be handled indirectly. \cite{kim2006} consider a rich VRP faced by a major waste management services provider in the United States. Although total cost minimization remains the main objective, equity issues concerning route compactness and workload balance are handled indirectly through a cluster-based solution procedure. The authors find that this approach outperforms classical insertion procedures when it comes to generating well-balanced routes. Similar observations are made in \cite{mendoza2009} in the context of periodic meter-maintenance routes for a public utility company in Colombia. In this case, workload is measured by the number of maintenance visits, and equity is evaluated based on the standard deviation of the number of visits per day. The authors find that allocating visits with a cluster-first, route-second approach leads to significantly better workload balance than generating routes first and allocating them to workdays afterwards. 

Finally, some applications use more sophisticated functions for measuring equity. \cite{jang2006} optimize the scheduling and routing of sales representatives of a lottery. Due to the central role of the representatives in acquiring and retaining customers, reducing overtime and increasing morale is a relevant concern. The authors measure workload equity with the mean square deviation of working hours and optimize a weighted sum of this measure combined with total travel distance. \cite{keskinturk2011} balance the tour durations for a baked goods distributor by optimizing an average relative imbalance measure as the primary objective. \cite{kritikos2010} consider the classic VRPTW with a weighted sum objective that includes a balanced cargo term, measured as the deviation from the median load. \cite{lee1999} measure equity as the sum of the differences to the minimum tour duration.

\section{Theoretical Analysis.}
\label{3-analysis}

In the preceding section, we observed that the VRPRB has been established in the literature as a prototypical model for an equitable VRP. The CVRP is one of the simplest VRPs, and the bi-objective model is a generalization of the weighted sum, constraint-based, and single-objective approaches. The VRPRB is thus a useful point of departure for further research.

Although many methods have already been proposed for the VRPRB, we could find no substantial discussion about the relative merits or limitations of the VRPRB model itself. Admittedly, the decision to use tour length as the equity \textit{metric} seems to be justified in practice, since most surveyed applications also base workload either fully or at least partly on tour length. However, the reasoning for and consequences of adopting the range as the equity \textit{function} in the VRPRB remain unclear, even though this is a defining feature of the model.

This raises the question to what extent the chosen objective affects the properties and structure of optimal VRP solutions. Some previous articles have examined this issue. \cite{bertazzi2015} consider the worst-case trade-off between the classic min-cost and min-max optimal solutions for various VRP variants. It is shown that in the worst case, the longest route of the optimal min-cost solution is \textit{k} times the longest route of the optimal min-max solution, and the min-max total distance is \textit{k} times the min-cost total distance, where~\textit{k} is the number of vehicles. \cite{campbell2008} study the effect of different objective functions in the context of humanitarian relief efforts, where the priority is to minimize the arrival time of aid packages. As reported by the authors, minimizing total routing cost leads to disproportionate increases in both the latest arrival time and the sum of arrival times. In contrast, minimizing the latest arrival time or the sum of arrival times leads to only moderate increases in total cost while significantly improving the service objectives. \cite{huang2012} examine weighted sum combinations of efficiency, efficacy, and equity objectives, and how they affect tour structure and resource utilization. These studies all conclude that the choice of objective function should be well-justified because it has important theoretical and practical implications.

The articles cited above clearly demonstrate the potential drawbacks of focusing on extreme single-objective optima, but worst-case analyses offer only very limited (and possibly too conservative) guidance with respect to the properties of intermediate trade-off solutions. Yet such insight is particularly relevant in practice. Indeed, even if cost minimization remains the main objective in many applications, more equitable solutions are sought and also accepted provided that the marginal cost of equity is reasonable. In such contexts, the primary interest lies ultimately in exploring various trade-off solutions rather than the extreme single-objective optima. A bi-objective analysis should therefore provide more generalizable conclusions.

To the best of our knowledge, \cite{halvorsen2016} present so far the only study examining the effects of different equity functions on VRP solutions in a bi-objective setting. The authors investigate a bi-objective MCGRP with a classical total cost objective and four possible equity objectives based on tour lengths: the range, the maximal tour, the total absolute deviation from the mean tour length, and the total absolute deviation to a target tour length. A worst-case analysis shows that if all tours are TSP-optimal, the total distance of the range-minimizing solution is up to \textit{k} times the min-cost total distance, where \textit{k} is the number of vehicles. Furthermore, computational experiments suggest that different measures lead to very different total costs even at maximal equity. When examining intermediate trade-off solutions, \cite{halvorsen2016} illustrate that there may be large differences in equity even among solutions with the same total cost, and vice versa. Finally, the results of their computational study suggest that the size of optimal Pareto sets may vary substantially depending on the chosen equity function.

Although a number of specific equity functions have been considered in all of the cited works, there is a lack of general guidance about what \textit{types} of functions may or may not be appropriate in a VRP context, and on what \textit{basis} to distinguish such type categories. Given the diversity of the surveyed applications, specific function recommendations can be of only limited scope, and thus general insights are needed. To close this gap, we complement and extend the existing works by approaching the choice of equity function from an axiomatic perspective, and deriving on this basis relevant properties of the corresponding Pareto optimal solution sets.

After briefly specifying the model assumptions in Section~\ref{assumptions}, we recall and discuss several widely-accepted axioms for equity measures in Section~\ref{axioms}, consider six classical measures in Section~\ref{measures}, and based on their properties, point out relevant implications for optimization in Section~\ref{implications}. To gauge the extent of certain implications and to assess other relevant aspects for choosing an equity measure, we report in Section~\ref{4-numerical} our observations from a numerical study on small instances solvable to optimality. We note that our conclusions are not always in agreement with previous works.


\subsection{Model Assumptions.}
\label{assumptions}

The VRPRB aims to capture the trade-off between efficiency/cost (total workload) and equity (the fair distribution of the workload). In this analysis, we assume that workload is to be minimized, and that each driver is allocated to exactly one tour (the tours in multi-trip VRPs may be suitably combined before allocating drivers). With respect to the evaluation of the two objectives, a VRP solution can be characterized as a vector $\mathbf{x}$ of $n$ workload allocations $x_1, \dotsc, x_n$, each corresponding to a different tour in the solution. In the VRPRB, the workload of a tour is given by its distance, so that the cost objective $C(\mathbf{x}) = \sum_{i=1}^{n} x_i$. What remains to be defined for the bi-objective problem is a function $I(\mathbf{x})$ which expresses the inequality of a solution's workload allocation.

It is clear that issues of equitable distribution are relevant only if the resource to be distributed is finite and the number of agents is at least two. Under these conditions, equitable distribution can be seen as a multi-objective problem with $n$ objectives, where each objective represents the allocation to a different agent. An equitable distribution should therefore be a Pareto-efficient solution to the $n$-objective allocation problem (\citealt{kostreva2004}). However, already for small $n$, the proportion of non-dominated and incomparable solutions rapidly increases, so that the Pareto-dominance relation is no longer helpful for distinguishing between different allocations (\citealt{farina2004}). Similarly, if the total workload to be distributed is constant, then every possible allocation is Pareto-efficient and incomparable. Thus, the need arises for measures which reduce the dimensionality of equity and introduce a stricter ordering among allocations..

\subsection{Desirable Properties of Inequality Measures.}
\label{axioms}

The choice of an inequality measure can be approached from an axiomatic perspective. The economics literature contains a large body of work on the subject of inequality measures and their desirable properties, mostly in the context of income distribution (e.g., ~\citealt{sen1973}, \citealt{allison1978}). Although there is no single accepted measure for assessing equity in general, a set of basic criteria have been identified which reasonable inequality measures should ideally satisfy. In the following we recall these criteria. We denote with $I(\mathbf{x})$ a function (inequality measure) which assigns an index value to any allocation $\mathbf{x}$ of $n$ outcomes (workloads) among the set of feasible allocations $X$. Without loss of generality, we assume that all outcomes $x_1, \dotsc, x_n$ as well as $I(\mathbf{x})$ itself are to be minimized.

\begin{axiom}
	(Inequality Relevance) If $x_i = x_j$ for all outcomes $i$ and $j$ in $\mathbf{x}$, then $I(\mathbf{x}) = 0$. Otherwise $I(\mathbf{x}) > 0$.
\end{axiom}

\begin{axiom}
	(Transitivity) Let $I(\mathbf{x}) \geq I(\mathbf{x'})$ and $I(\mathbf{x'}) \geq I(\mathbf{x''})$, then $I(\mathbf{x}) \geq I(\mathbf{x}'')$.
\end{axiom}

The criterion of inequality relevance simply states that the measure should have a value of 0 (or any other known and fixed value) if the distribution is perfectly equitable, and be positive otherwise. Clearly, a measure should be transitive so that the ordering it produces is consistent with itself. Most inequality measures satisfy these two criteria.

\begin{axiom}
	(Scale Invariance) $I(\mathbf{x}) = I(\lambda \mathbf{x})$ for $\lambda \in \mathbb{R}\setminus\{0\}$.
\end{axiom}

\begin{axiom}
	(Translation Invariance) Let $\alpha \in \mathbb{R}$, and $\mathbf{u}$ be a unit vector of length $n$. Then $I(\mathbf{x}) = I(\mathbf{x} + \alpha \mathbf{u})$.
\end{axiom}

Scale invariance requires that the inequality index remains unchanged if all outcomes in a distribution are multiplied by a constant. This ensures that the degree of inequality of an allocation is not affected by the unit of measure (e.g.~measuring distance in kilometres or miles). In contrast, translation invariance requires that the inequality index remains unchanged if a constant is added to every outcome in the distribution. Both of these properties reflect specific value judgements.

\begin{axiom}
	(Population Independence) $I(\mathbf{x}) = I(\mathbf{x} \cup \mathbf{x} \cup \dotsc \cup \mathbf{x})$.
\end{axiom}

Some measures are affected by the number of outcomes in a distribution. In such cases, a direct comparison of distributions of different sizes is not possible. One way to circumvent this limitation is to replicate each population a certain number of times such that both resulting populations are of the same size (e.g.~the lowest common multiple of the original sizes). Such transformations are not necessary if the measure is population independent \citep{aboudi2010}.

\begin{axiom}
	(Anonymity) Let $\mathbf{x'}$ be a permutation of the elements in $\mathbf{x}$, then $I(\mathbf{x}) = I(\mathbf{x'})$.
\end{axiom}

Also known as \textit{impartiality} or \textit{symmetry}, this axiom states that a distribution's equity should not depend on a particular ordering of the outcomes. This assumption may seem self-evident, but if allocation is performed over entities which are not identical (different needs, rights, preferences), then an anonymous inequality measure is not appropriate without some additional treatment. This can be accomplished by introducing what \cite{marsh1994} refer to in their framework as a \textit{reference distribution} specifying the target outcome for each entity. Equity then depends on the individual deviations to the respective targets, and if these targets differ, then a perfectly equal distribution of outcomes is not optimal. In this context, the survey by \cite{karsu2015} distinguishes between equity (anonymity holds) and balancing (anonymity does not hold), and includes a review of optimization models in which anonymity does not hold. Practically all the surveyed VRP publications deal with anonymous models.

\begin{axiom}
(Monotonicity) Let $\mathbf{x'}$ be formed as follows: $x'_i = x_i + \delta_i$ \small{for at least one} $i$ in $\mathbf{x}$. If $\delta_i \geq 0$ for all $i$ with at least one strict inequality, then $I(\mathbf{x'}) \geq I(\mathbf{x})$.
\end{axiom}

Recall that equitable distribution among $n$ recipients may be seen as an $n$-objective optimization problem. It thus follows that the inequality measure should prefer \textit{Pareto-optimal} distributions \citep{marsh1994, eiselt1995}, or in other words, be consistent also with the optimization of each individual outcome \citep{ogryczak2014}. Formally, this means that the inequality measure should be at least weakly monotonic in all outcomes of the distribution.

\begin{axiom}
 (Pigou-Dalton Transfer Principle) Let $\mathbf{x'}$ be formed as follows: $x'_i = x_i + \delta$, $x'_j = x_j - \delta$, $x'_k = x_k$, for all $k \notin \{i,j\}$. If $0 \leq \delta < x_j - x_i$, then $I(\mathbf{x'}) \leq I(\mathbf{x})$.
\end{axiom}

The Pigou-Dalton principle of transfers (PD) is a widely-accepted criterion for inequality measures when the aforementioned anonymity assumption holds. A transfer refers to shifting $\delta$ units of workload from an individual $j$ to another individual $i$, such that $x'_i = x_i + \delta$ and $x'_j = x_j - \delta$ (by definition, such transfers are \textit{mean-preserving}). A transfer is said to be \textit{progressive} (favouring the worse-off party) if $0 \leq \delta \leq x_j - x_i$, and \textit{regressive} otherwise. The weak version of the PD principle states that if an allocation $\mathbf{x'}$ can be reached by a finite series of only progressive transfers from an allocation $\mathbf{x}$, then $I(\mathbf{x'}) \leq I(\mathbf{x})$ (i.e.~the new allocation cannot be less equitable). The strong version of the PD principle requires strict inequality, i.e.~$I(\mathbf{x'}) < I(\mathbf{x})$. Note that the PD principle applies only if the allocations do not differ in the number and the sum of their outcomes.



\subsection{Commonly Used Inequality Measures.}
\label{measures}

Based on our survey, the review by \cite{karsu2015}, and the broader economics literature, we have identified six classical measures of equity. In order of sophistication, they are: the worst outcome (in our context minimization of the maximal workload, which we will refer to as min-max), the lexicographic extension of min-max, the range, the mean absolute deviation, the standard deviation, and the Gini coefficient. In the following we briefly discuss these measures with respect to the axioms introduced above. A summary overview is provided in Table~\ref{measuresOverview}.

\begin{table}[h]
\renewcommand{\arraystretch}{1.2}
\setlength{\tabcolsep}{3.5pt}
\caption{Common Inequality Measures and their Axiomatic Properties}
\hspace*{+0.25cm}
\scalebox{0.9}
{
\begin{tabular}{|l|M{2.1cm}|M{2.1cm}|M{2.1cm}|M{2.1cm}|M{2.1cm}|c|}\hline
		\textbf{Property} 					& \textbf{Min-max} 	& \textbf{Lexicogr. min-max} 	& \textbf{Range} & \textbf{Mean abs. deviation} & \textbf{Standard deviation} 	& \textbf{Gini coeff.} 	\\ \hline
		Inequality Rel. 	&   				& 	 					& $\bullet$ & $\bullet$	& $\bullet$ & $\bullet$					\\ 
    Transitive 				& $\bullet$ & $\bullet$		& $\bullet$	& $\bullet$ & $\bullet$ & $\bullet$					\\ 
    Scale Inv. 				&   				&   					&   				&   				&   				& $\bullet$  				\\ 
		Translation Inv. 	&   				&   					& $\bullet$ & $\bullet$ & $\bullet$ &    								\\ 
		Population Ind. 	& $\bullet$	& $\bullet$ 	& $\bullet$ & $\bullet$	& $\bullet$ & $\bullet$  				\\ 
		Anonymous 				& $\bullet$ & $\bullet$ 	& $\bullet$ & $\bullet$ & $\bullet$ & $\bullet$  				\\ 
		Monotonicity	 		& weak	 		& strong			& 		  		& 				 	& 			 		& 			    				\\
		Pigou-Dalton	 		& weak	 		& strong			& weak 			& weak(+)		& strong 		& strong    				\\  \hline
\end{tabular}
}
\label{measuresOverview}
\end{table}

\paragraph{\bf Min-Max:}

$\max x_i$

Optimizing the worst outcome is likely the simplest inequality measure there is. In addition to the min-max VRP, it is also the main objective in the class of \textit{p}-center facility location problems as well as makespan minimization problems in scheduling (though equity concerns are technically not the focus in the latter). A justification of using min-max as an equity measure may be found in \cite{rawls}, in the context of distributive justice. If an allocation problem is anonymous and impartial, then none of the entities subject to the allocation know their position in the distribution. From this point of view (termed \textit{veil of ignorance}), individuals would rather prefer min-max distributions as they remain acceptable even in the worst case.

However, the simplicity of the min-max approach has clear disadvantages. It cannot distinguish between distributions with identical worst outcomes (min-max would be indifferent between [20,15,10,5] and [20,10,10,10], though the second is clearly more equitable), and if the worst outcomes do differ, min-max can neither capture nor quantify the differences in equitability of the remaining parts of the distribution (e.g.~[19,19,11,1] would be considered preferable to [20,10,10,10]). Note that this remains true even if the total distributed workload is constant, as in the previous examples.

\vspace*{0.3cm}
\paragraph{\bf Lexicographic Min-Max:}

$ $

The lexicographic approach is a natural extension of the min-max principle. Rather than minimizing only the worst outcome, also the second-worst is minimized (subject to minimization of the first), the third-worst is minimized (subject to minimization of the previous two), and so on. The lexicographic measure satisfies the strong version of the PD principle and provides a strict preference ordering such that every distribution is assigned a unique rank in the order. This resolves the problem of distinguishing between distributions with identical worst outcomes, but does not quantify the remaining differences. Indeed, the lexicographic rank alone does not permit to make statements about \textit{how much more} equitable one distribution is compared to another.

\vspace*{0.3cm}
\paragraph{\bf Range:}

$\max x_i - \min x_i$

The range is the difference between the worst and best outcomes. Despite being the simplest measure of dispersion, it already provides more information than the simple min-max. At the same time, it remains simple to implement and interpret, making it a popular choice in applications. Another desirable property is that the range has an unambiguous optimal value of 0, unlike min-max and its lexicographic version. However, the range does not capture the absolute levels of the outcomes (e.g.~[10,9,8,7] would be preferred over [5,4,2,1], even though all workloads are lower in the latter). In addition, the range satisfies only the weak version of the PD principle, and any changes between the extremes have no effect on the index. This motivates the use of more sophisticated measures.

\vspace*{0.3cm}
\paragraph{\bf Mean Absolute Deviation:}

$\frac{1}{n}\sum_{i=1}^{n} |x_i-\bar{x}|$

The mean absolute deviation (MAD) is defined as the mean absolute difference between each outcome and the mean outcome. Unlike the range, MAD is directly affected by every outcome in the distribution, rather than the extremes only. However, it does not satisfy the strong version of the PD principle, because transfers between outcomes on the same side of the mean do not affect this measure. Nonetheless, this is a smaller proportion of transfers than all those which do not affect the range, so in that sense, MAD is 'stronger' than the range with respect to the PD principle.

\vspace*{0.3cm}
\paragraph{\bf Standard Deviation:}

$\sqrt{\frac{1}{n} \sum_{i=1}^n (x_i - \bar{x})^2}$

The standard deviation is arguably the most well-known statistical measure of dispersion. It satisfies nearly all the aforementioned axioms. In particular, it satisfies the strong version of the PD principle, unlike three of the previous four measures. The standard deviation is translation invariant, but can easily be transformed into a scale invariant measure by dividing it by the mean of the distribution, which yields the \textit{coefficient of variation}. Although it does not provide a strict ordering of distributions, it can be expected that in practice, distributions with identical standard deviations will be rare. The main disadvantage of the standard deviation is its computational complexity, and for decision makers it can be a less intuitive measure than simpler ones such as the range.

\vspace*{0.3cm}
\paragraph{\bf Gini Coefficient:}

$\frac{1}{2n^2\bar{x}} \sum_{i=1}^{n} \sum_{j=1}^{n} |x_i-x_j|$

The Gini coefficient is one of the most widely used indices in income economics and for inequality studies in general. The index assumes values between 0 and 1, and is easy to interpret: lower values correspond to lower inequality. In practice, the upper bound depends on the number of outcomes and is therefore less than 1 for a finite set of outcomes. A translation invariant version of the index can be obtained by multiplying it with the mean of the outcomes. Like the standard deviation, the Gini coefficient satisfies the strong version of the PD principle, and similarly, the complexity of this measure is its main disadvantage in optimization settings. However, it is worth noting that the Gini index can be formulated as a linear objective \citep{mandell1991}, and that it becomes linear in the special case when the average outcome $\bar{x}$ and the number of outcomes $n$ are both constant.

\subsection{Implications for Choosing an Equity Objective.}
\label{implications}

Based on the foregoing analysis, we can see that none of the examined measures comply with every axiom, and no single measure dominates all others in every aspect. Depending on the monotonicity property, they can be divided into two groups: \textit{absolute} measures (min-max and lexicographic min-max), and \textit{relative} measures (the remaining four). Both groups include at least one measure that satisfies the strong version of the PD principle.

With respect to the PD principle, it is useful to distinguish between two general types of metrics: those whose sum is constant over all solutions (e.g.~demand or number of customers in typical VRPs), and those for which the sum is variable (e.g.~tour lengths).

\vspace*{0.3cm}
\begin{definition} {\bf Constant / Variable Sum Equity Metric:}
\label{def1}
An equity metric is constant-sum if $\sum_{i=1}^{n} x_i$ is identical for all solutions $\mathbf{x} \in X$, and variable-sum otherwise.
\end{definition}
\vspace*{0.3cm}

By definition, the PD principle applies only when comparing solutions with the same sum of outcomes. This leads directly to the following observation:

\vspace*{0.2cm}
\begin{observation}
\label{R1}
The Pigou-Dalton transfer principle does not, in general, affect preferences between solutions if the equity metric is variable-sum.
\end{observation}
\vspace*{0.2cm}

The extent to which this holds depends on the particular problem and instance, and specifically on the prevalence of distinct solutions with identical values for one objective but not for the other. Of course, the \textit{opposite} observation can be made when equity is based on constant-sum metrics. In such cases, every pair of feasible solutions is connected by mean-preserving transfers, and differences in equity can be captured with the PD principle whenever those transfers are strictly regressive. Since tour lengths are variable-sum, an inequality measure's satisfaction or not of the PD principle is not expected to have any noticeable impact on the resulting efficient sets for the standard VRPRB.

On the other hand, the monotonicity property has several direct consequences for the VRPRB if workload is modelled using tour lengths or some other variable-sum metric, as done in most of the surveyed papers. We identify the following implications:

\vspace*{0.2cm}
\begin{theorem}
\label{T1}
If $I(\mathbf{x})$ is at least weakly monotonic, then every Pareto-optimal solution to the VRPRB is composed only of TSP-optimal tours.
\end{theorem}
\vspace*{0.2cm}

\proof{Proof.}
By contradiction. Let $\mathbf{x}$ be a Pareto-optimal solution to the VRPRB which contains at least one non-TSP-optimal tour, then there exists a solution $\mathbf{x'}$  composed of the corresponding TSP-optima of all tours in $\mathbf{x}$. This implies $C(\mathbf{x}) > C(\mathbf{x'})$ and $I(\mathbf{x}) \geq I(\mathbf{x'})$ since $I(\mathbf{x})$ is at least weakly monotonic in all $x_i$, but this contradicts the Pareto-optimality of $\mathbf{x}$.~\Halmos
\endproof

In contrast, no such guarantee exists for non-monotonic measures. By definition, non-monotonicity implies that the deterioration of an outcome (e.g.~non-TSP-optimality of a tour) can worsen, leave unchanged, but also improve the value returned by the inequality index. The latter case represents a potentially Pareto-optimal solution increasing cost while seemingly improving equity.

\vspace*{0.2cm}
\begin{observation}
\label{T2}
If $I(\mathbf{x})$ is not monotonic, then the Pareto-optimal solutions of the VRPRB can include non-TSP-optimal tours.
\end{observation}
\vspace*{0.2cm}


Clearly this is not a desirable property, as non-TSP-optimal tours contradict the minimization of individual workloads. More generally, non-monotonic objectives can be improved even in the extreme case when \textit{all} workloads are strictly increased. This issue has previously been observed in the field of facility location, where it is known that the direct optimization of non-monotonic measures may lead to locating the facilities as far away as possible, providing almost perfectly equal \textit{lack} of service \citep{marsh1994, eiselt1995, ogryczak2000}. 

Contrary to intuition, adding an explicit total cost objective does not resolve this issue, because an improvement in even just one objective (in this case equity) is already a sufficient condition for Pareto-optimality. We introduce the concept of \textbf{(workload) inconsistency} to distinguish solutions of this type.

\vspace*{0.2cm}
\begin{definition} {\bf Inconsistency:}
\label{def2}
A solution $\mathbf{x}$ is inconsistent if there exists another solution $\mathbf{x'} \in X$ such that $C(\mathbf{x}) > C(\mathbf{x'})$ and $x_i \geq x'_i$ for all $i \in \{1, \dots, n\}$, but $I(\mathbf{x}) < I(\mathbf{x'})$.
\end{definition}
\vspace*{0.2cm}

It follows that lower bound constraints on workload (as proposed in some papers) can directly lead to inconsistency, regardless of the equity function. Perhaps more surprisingly, even TSP-optimal solutions can be inconsistent, since a different, possibly very-sub-optimal allocation of customers to tours can also increase all tour lengths. This leads to the following observation:

\vspace*{0.2cm}
\begin{observation}
\label{P1}
TSP-optimality is not a sufficient condition for avoiding inconsistency.
\end{observation}
\vspace*{0.2cm}

A relevant question that remains is how monotonicity of the equity function relates to the inconsistency of the generated Pareto-optimal solutions.

\vspace*{0.2cm}
\begin{theorem}
\label{T3}
If $I(\mathbf{x})$ is at least weakly monotonic, then the Pareto-optimal solutions of the VRPRB cannot be inconsistent.
\end{theorem}
\vspace*{0.2cm}

\proof{Proof.}
By contradiction. Let $\mathbf{x}$ be an inconsistent Pareto-optimal solution to the VRPRB. According to Definition~\ref{def1}, there must exist another solution $\mathbf{x'} \in X$ such that $C(\mathbf{x}) > C(\mathbf{x'})$ and $x_i \geq x'_i$ for all $i \in \{1, \dots, n\}$, but $I(\mathbf{x}) < I(\mathbf{x'})$. If $I(\mathbf{x})$ is at least weakly monotonic and $x_i \geq x'_i$ for all $i \in \{1, \dots, n\}$, then this implies $I(\mathbf{x}) \geq I(\mathbf{x'})$, which contradicts the previous statement and Definition~\ref{def1}.~\Halmos
\endproof

\begin{observation}
\label{T4}
If $I(\mathbf{x})$ is not monotonic, then the Pareto-optimal solutions of the VRPRB can be inconsistent.
\end{observation}


Combining the results of Theorems~\ref{T1} and~\ref{T3}, it follows that monotonic measures of equity guarantee both TSP-optimality and consistency of all generated Pareto-optimal solutions. Neither of these properties can be guaranteed for non-monotonic measures without additional restrictions. With non-monotonic inequality measures, it is therefore \textit{possible} that some Pareto-optimal VRPRB solutions will be non-TSP-optimal and inconsistent. However, it remains unclear to what \textit{extent} inconsistency and non-TSP-optimality are a relevant concern in practice for non-monotonic measures.

There are also a number of other aspects that are relevant for the choice of an equity measure, but which are not directly implied by its axiomatic properties. A practical issue may be the extent to which different measures generate the same Pareto-optimal solutions. If there is considerable agreement between some measures, then one may prefer the simpler one for reasons of computational complexity and ease of understanding. Likewise, if there is a noticeable difference in the average cardinality of the efficient sets, then this has practical implications (the preference for smaller or larger fronts will also depend on the application). Finally, we are also interested in the marginal cost of equity, i.e.~the extent to which cost minimization and equity are conflicting objectives.

In the next section, we provide some insight into these aspects by conducting a computational study on a set of instances solvable to optimality.
\section{Numerical Analysis.}
\label{4-numerical}

In this section, we assess numerically how the choice of equity function affects the resulting Pareto sets and the structure of Pareto-optimal VRPRB trade-off solutions. Our experiments are based on a set of small CVRP instances solvable to optimality for all of the examined measures.

We derived our test instances from the benchmark instance CMT5, which contains all the customer coordinates of the instances CMT1 to CMT10 of \cite{CMT}. To the best of our knowledge, the customers are distributed randomly and without any special structure. We used the first $n$ customers to generate a subset of instances by varying the number of vehicles $v$ between 2 and 5, and setting the capacity $q$ to either $\left\lceil \frac{n}{v} \right\rceil$ (constrained) or $\left\lceil \frac{n}{v} \right\rceil + 1$ (less constrained). Unit capacity was used for all customers. The next subset of instances was created using the next $n$ customers of CMT5, and so on. We generated the optimal Pareto sets by computing all feasible tours and then examining all feasible combinations of tours to derive all feasible solutions. Our study is based on instances with $n=14$, as these were the largest for which Pareto-optimal sets for all parameter combinations could be enumerated within a reasonable time. We base our analysis on a total of 60 instances falling into 10 disjoint customer subsets. 

\subsection{General Observations.}

\vspace*{0.2cm}
\paragraph{\bf Equity of Cost-Optimal Solutions.}
Previous studies have noted that cost-optimal VRP solutions tend to be poorly balanced \citep{jozefowiez2009,reiter2012,bertazzi2015}. We confirm this and find that on average, the longest tour is about twice as long as the shortest. This motivates the search for alternative trade-off solutions. 

\vspace*{0.2cm}
\paragraph{\bf Marginal Cost of Equity.}
We observe that the marginal trade-off between cost and equity is generally favourable: by selecting the TSP-optimal and consistent solution with second-best cost, the difference between the longest and shortest tours can already be reduced by around 40\% on average, while total cost increases by only around 2\%. In addition, about 37\% of efficient solutions were found to have total costs within 10\% of the cost optimum. These observations suggest that reasonable compromise solutions can be found in practice.


\vspace*{0.2cm}
\paragraph{\bf Agreement Between Equity Functions.} On average, the highest degree of agreement, slightly over 80\%, occurs unsurprisingly between the min-max measure and its lexicographic extension, since the solutions generated by the former are necessarily a subset of the latter. The degree of overlap among the non-monotonic measures is lower, from around 50\% between the range and the mean absolute deviation, to around 80\% between the standard deviation and Gini coefficient. These are rather moderate degrees of agreement considering the size of the instances, and this agreement is likely to decline significantly with larger instance size. Overall, only a small subset of solutions (around 14\% in our study) is identified by \textit{all} examined measures. In contrast, a quarter of solutions is unique to no more than two measures: around 14\% are unique to only one, and around 10\% to only two. Choosing which measure to use is thus a meaningful decision.

\begin{figure}[H]
\FIGURE
{\includegraphics[width=1.0\textwidth]{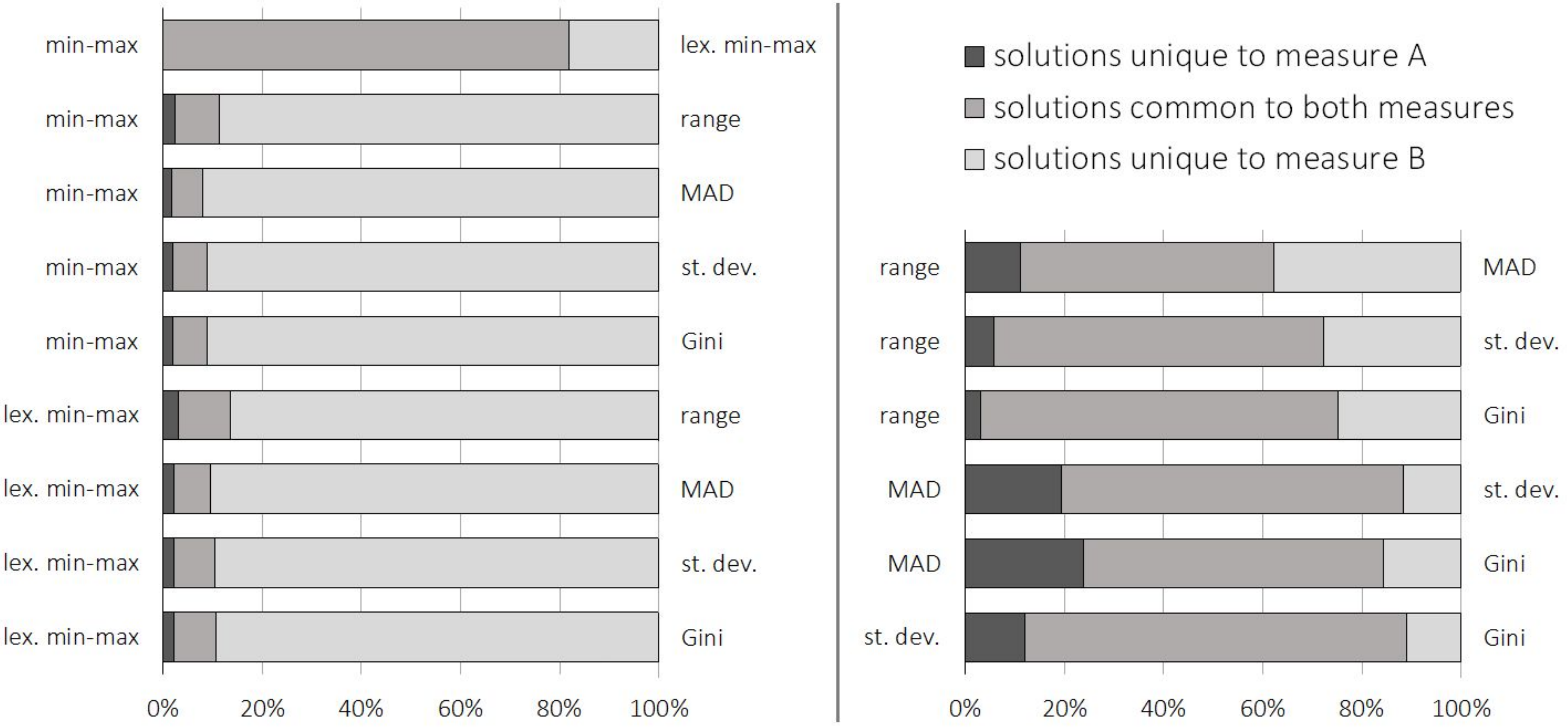}}
{Agreement Between the Examined Equity Functions.
\label{overlap}}
{}
\end{figure}

We observe that a conventional solution of the instances produces, at least for non-monotonic equity functions, optimal Pareto sets which contain solutions with non-TSP-optimal tours. The VRPRB literature proposes to apply a 2-opt local search to each solution prior to examining Pareto efficiency. In order to assess the extent of this phenomenon and the impact of adding extra constraints to the model, we solve all instances twice: first conventionally, and then with the constraint that all tours must be TSP-optimal. 

\subsection{Conventional Optimization.}

We first solve each VRPRB instance to optimality without any further restriction (e.g.~2-optimality) on the acceptance of solutions. In Table~\ref{summary1}, we report for each equity function the average cardinality of the optimal Pareto sets, the average number of TSP-optimal solutions in those sets, and the average number of consistent solutions, in absolute terms and relative to the average total. In order to give a more detailed impression of the differences between individual equity functions, we present in Table~\ref{sample1tabular} the complete Pareto fronts of all examined equity functions for a representative instance. Table~\ref{sample1tabular} lists for each solution the total cost, the solution's tour lengths in decreasing order, for each examined equity function the solution's corresponding inequality value, whether the solution is TSP-optimal, and whether the solution is consistent (white circle) or not (black circle) in the corresponding equity function's optimal front. A dash denotes that the solution is not Pareto-optimal for the respective equity function and therefore absent from that function's optimal Pareto set.

\begin{table}[htbp]
\centering
\SingleSpacedXI
\caption{Average Number of Pareto Efficient, TSP-optimal, and Consistent Solutions per Equity Function -- Conventional Optimization}
\renewcommand{\arraystretch}{1.2}
\resizebox{\textwidth}{!}{%
\begin{tabular}{lrcrcrcrcrcrc} \hline
                      & \multicolumn{2}{M{2.5cm}}{\textbf{min-max}} & \multicolumn{2}{M{2.5cm}}{\textbf{lexicogr. min-max}} & \multicolumn{2}{M{2.5cm}}{\textbf{range}} & \multicolumn{2}{M{2.5cm}}{\textbf{mean abs. deviation}} & \multicolumn{2}{M{2.5cm}}{\textbf{standard deviation}} & \multicolumn{2}{M{2.5cm}}{\textbf{Gini coefficient}} \\ \hline
                      & abs.	& \%	& abs.	& \%	& abs.	& \%	& abs.	& \%	& abs.	& \%	& abs.	& \%	\\ \hline
\textbf{Total}				& 3.53	& -			& 4.23	& - 		& 30.22	& -			& 43.18	& -			& 39.37	& -			& 38.97	& -  		\\
\textbf{TSP-optimal}  & 3.53	& 100\% 	& 4.23 	& 100\%		& 5.58	& 18\% 		& 6.62	& 	15\%	& 6.20 	& 	16\%  & 6.03	& 	15\%	\\
\textbf{Consistent}		& 3.53	& 100\%		& 4.23	& 100\%		& 4.32	& 14\%		& 4.55	& 	11\%	& 4.52	& 	11\% 	& 4.65	& 	12\%	\\ \hline
\end{tabular}
}
\label{summary1}
\end{table}

Based on the summary statistics in Table~\ref{summary1}, clear distinctions in behaviour can be observed between the monotonic and non-monotonic measures. There is a marked difference in the average cardinality of the optimal Pareto sets. The non-monotonic measures generate around 10 times as many Pareto-efficient solutions as the monotonic ones. In fact, the resulting Pareto fronts are unexpectedly large (around 30 to 40 solutions) considering the size of the instances. However, we also observe that only a fairly small percentage (15\% to 18\%) of these efficient solutions are TSP-optimal and even fewer are consistent (11\% to 14\%), in contrast to 100\% with the monotonic measures. This is in agreement with the implications identified in Section \ref{implications} and shows that the extent of non-TSP-optimality and inconsistency should not be underestimated.

\begin{landscape}
\begin{table}[htbp]
\SingleSpacedXI
\renewcommand{\arraystretch}{1.2}
\vspace*{0.5cm}
\caption{Pareto Optimal Fronts per Equity Function for a Representative Instance (\#30), Conventional Optimization}
\label{sample1tabular}
\hspace*{-1.2cm}
\scalebox{0.85}
{
\begin{tabular}{|cF|ccccc|FFFFFF|c|cccccc}
\hline
 & \textbf{} &  \multicolumn{5}{c|}{\textbf{Tour Lengths}} & \multicolumn{6}{c|}{\textbf{Equity Function}} & \textbf{}  & \multicolumn{6}{c|}{\textbf{Inconsistency}} \\ \hline

\multicolumn{1}{|c|}{\textbf{Nr.}} & \textbf{Total Cost} & \textbf{Tour 1} & \textbf{Tour 2} & \textbf{Tour 3} & \textbf{Tour 4} & \textbf{Tour 5} & \textbf{minmax (E1)} & \textbf{lex.~rank (E2)} & \textbf{range (E3)} & \textbf{MAD (E4)} & \textbf{std.~dev. (E5)} & \textbf{Gini (E6)} &  \textbf{TSP} & \textbf{E1} & \textbf{E2} & \textbf{E3} & \textbf{E4} & \textbf{E5} & \multicolumn{1}{c|}{\textbf{E6}} \\ \hline
\multicolumn{1}{|c|}{\textbf{1}}	 & 340.45              & 97.15           & 78.04           & 68.96           & 53.13           & 43.17           & 97.15           & 1              	  & 53.98          & 15.95                & 18.92              & 0.16          & $\checkmark$         & $\circ$     & $\circ$     & $\circ$     & $\circ$     & $\circ$     & \multicolumn{1}{c|}{$\circ$}            \\
\multicolumn{1}{|c|}{\textbf{2}}   & 341.83              & 83.04           & 78.04           & 68.96           & 68.61           & 43.17           & 83.04           & 2               & 39.87          & 10.08                & 13.74              & 0.10          & $\checkmark$         & $\circ$     & $\circ$     & $\circ$     & $\circ$     & $\circ$     & \multicolumn{1}{c|}{$\circ$}            \\
\multicolumn{1}{|c|}{\textbf{3}}   & 344.31              & 83.04           & 78.04           & 71.44           & 68.61        	 & 43.17           & -               & -                  & -              & -                    & -                  & 0.10          &              	&             &             &             &             &             & \multicolumn{1}{c|}{$\bullet$}           \\
\multicolumn{1}{|c|}{\textbf{4}}   & 345.43              & 83.04           & 77.81           & 72.80           & 68.61           & 43.17           & -               & 3               & -              & -                    & -                  & 0.10          & $\checkmark$         &             & $\circ$     &             &             &             & \multicolumn{1}{c|}{$\circ$}            \\
\multicolumn{1}{|c|}{\textbf{5}}   & 346.41              & 83.04           & 78.04           & 73.20           & 68.96        	 & 43.17           & -               & -                  & -              & -                    & -                  & 0.10          &              		&             &             &             &             &             & \multicolumn{1}{c|}{$\bullet$}           \\
\multicolumn{1}{|c|}{\textbf{6}}   & 347.61              & 83.04           & 75.52           & 72.80           & 68.61           & 47.64           & -               & 4               & 35.40          & 9.12                 & 11.91              & 0.09          & $\checkmark$         &             & $\circ$     & $\circ$     & $\circ$     & $\circ$     & \multicolumn{1}{c|}{$\circ$}            \\
\multicolumn{1}{|c|}{\textbf{7}}   & 350.44              & 83.04           & 75.48           & 72.80           & 68.61           & 50.50           & -               & 5               & 32.54          & 8.42                 & 10.86              & 0.08          & $\checkmark$         &             & $\circ$     & $\circ$     & $\circ$     & $\circ$     & \multicolumn{1}{c|}{$\circ$}            \\
\multicolumn{1}{|c|}{\textbf{8}}   & 352.67              & 83.04           & 68.96           & 68.61           & 66.37           & 65.68           & -               & 6               & 17.36          & 5.00                 & 6.38               & 0.04          & $\checkmark$         &             & $\circ$     & $\circ$     & $\circ$     & $\circ$     & \multicolumn{1}{c|}{$\circ$}            \\
\multicolumn{1}{|c|}{\textbf{9}}   & 354.26              & 82.81           & 81.28           & 78.04           & 68.96           & 43.17           & 82.81           & 7               & -              & -                    & -                  & -             & $\checkmark$         & $\circ$     & $\circ$     &             &             &             & \multicolumn{1}{c|}{}            \\
\multicolumn{1}{|c|}{\textbf{10}}  & 355.15              & 83.04           & 71.44           & 68.61           & 66.37        	 & 65.68           & -               & -                  & -              & 4.97                 & 6.33               & -             &              		&             &             &             & $\bullet$   & $\bullet$   & \multicolumn{1}{c|}{}            \\
\multicolumn{1}{|c|}{\textbf{11}}  & 355.89              & 83.04           & 69.63           & 68.96           & 68.61        	 & 65.65           & -               & -                  & -              & 4.75                 & 6.09               & 0.04          & $\checkmark$         &             &             &             & $\circ$     & $\circ$     & \multicolumn{1}{c|}{$\circ$}            \\
\multicolumn{1}{|c|}{\textbf{12}}  & 356.07              & 81.28           & 78.04           & 74.67           & 68.96           & 53.13           & 81.28           & 8               & -              & -                    & -                  & -             & $\checkmark$         & $\circ$     & $\circ$     &             &             &             & \multicolumn{1}{c|}{}            \\
\multicolumn{1}{|c|}{\textbf{13}}  & 357.45              & 78.04           & 74.67           & 68.96           & 68.61           & 67.16           & 78.04           & 9               & 10.88          & 3.89                 & 4.16               & 0.03          & $\checkmark$         & $\circ$     & $\circ$     & $\circ$     & $\circ$     & $\circ$     & \multicolumn{1}{c|}{$\circ$}            \\
\multicolumn{1}{|c|}{\textbf{14}}  & 359.93              & 78.04           & 74.67           & 71.44           & 68.61           & 67.16        	 & -               & -                  & -              & 3.50                 & 3.97               & 0.03          &                 &             &             &             & $\bullet$   & $\bullet$   & \multicolumn{1}{c|}{$\bullet$}           \\
\multicolumn{1}{|c|}{\textbf{15}}  & 361.05              & 77.81           & 74.67           & 72.80           & 68.61           & 67.16           & 77.81           & 10              & 10.65          & 3.46                 & 3.90               & 0.03          & $\checkmark$         & $\circ$     & $\circ$     & $\circ$     & $\circ$     & $\circ$     & \multicolumn{1}{c|}{$\circ$}            \\
\multicolumn{1}{|c|}{\textbf{16}}  & 362.47              & 78.04           & 74.67           & 72.19           & 68.96           & 68.61        	 & -               & -                  & 9.44           & 3.09                 & 3.56               & 0.03          & $\checkmark$         &             &             & $\bullet$   & $\bullet$   & $\bullet$   & \multicolumn{1}{c|}{$\bullet$}           \\
\multicolumn{1}{|c|}{\textbf{17}}  & 364.23              & 78.04           & 74.67           & 72.19           & 70.36           & 68.96        	 & -               & -                  & 9.09           & 2.81                 & 3.23               & 0.03          &                 &             &             & $\bullet$   & $\bullet$   & $\bullet$   & \multicolumn{1}{c|}{$\bullet$}           \\
\multicolumn{1}{|c|}{\textbf{18}}  & 364.95              & 78.04           & 74.67           & 72.19           & 71.44           & 68.61        	 & -               & -                  & -              & 2.69                 & 3.18               & 0.02          &                 &             &             &             & $\bullet$   & $\bullet$   & \multicolumn{1}{c|}{$\bullet$}           \\
\multicolumn{1}{|c|}{\textbf{19}}  & 365.64              & 77.81           & 74.67           & 73.20           & 72.80           & 67.16        	 & -               & -                  & -              & 2.52                 & -                  & -             &                 &             &             &             & $\bullet$   &             & \multicolumn{1}{c|}{}            \\
\multicolumn{1}{|c|}{\textbf{20}}  & 366.08              & 77.81           & 74.67           & 72.80           & 72.19           & 68.61        	 & -               & -                  & -              & 2.42                 & 3.02               & 0.02          & $\checkmark$         &             &             &             & $\bullet$   & $\bullet$   & \multicolumn{1}{c|}{$\bullet$}           \\
\multicolumn{1}{|c|}{\textbf{21}}  & 366.71              & 78.04           & 74.67           & 72.19           & 71.44           & 70.36        	 & -               & -                  & 7.68           & 2.41                 & 2.75               & 0.02          &                 &             &             & $\bullet$   & $\bullet$   & $\bullet$   & \multicolumn{1}{c|}{$\bullet$}           \\
\multicolumn{1}{|c|}{\textbf{22}}  & 367.83              & 77.81           & 74.67           & 72.80           & 72.19           & 70.36        	 & -               & -                  & 7.45           & 2.14                 & 2.53               & 0.02          &                 &             &             & $\bullet$   & $\bullet$   & $\bullet$   & \multicolumn{1}{c|}{$\bullet$}           \\
\multicolumn{1}{|c|}{\textbf{23}}  & 368.96              & 77.40           & 75.48           & 74.67           & 72.80           & 68.61           & 77.40           & 11              & -              & -                    & -                  & -             & $\checkmark$         & $\circ$     & $\circ$     &             &             &             & \multicolumn{1}{c|}{}            \\
\multicolumn{1}{|c|}{\textbf{24}}  & 369.76              & 75.48           & 74.67           & 74.60           & 72.80           & 72.21           & 75.48           & 12              & 3.27           & 1.16                 & 1.24               & 0.01          & $\checkmark$         & $\circ$     & $\circ$     & $\circ$     & $\circ$     & $\circ$     & \multicolumn{1}{c|}{$\circ$}            \\
\multicolumn{1}{|c|}{\textbf{25}}  & 372.01              & 75.48           & 74.67           & 74.60           & 74.46           & 72.80        	 & -               & -                  & 2.69           & 0.64                 & 0.88               & 0.01          &                 &             &             & $\bullet$   & $\bullet$   & $\bullet$   & \multicolumn{1}{c|}{$\bullet$}           \\
\multicolumn{1}{|c|}{\textbf{26}}  & 405.72              & 81.90           & 81.45           & 81.32           & 81.28           & 79.79        	 & -               & -                  & 2.11           & 0.54                 & 0.71               & 0.00          &                 &             &             & $\bullet$   & $\bullet$   & $\bullet$   & \multicolumn{1}{c|}{$\bullet$}           \\
\multicolumn{1}{|c|}{\textbf{27}}  & 407.69              & 81.90           & 81.75           & 81.45           & 81.32           & 81.28        	 & -               & -                  & 0.62           & 0.23                 & 0.25               & 0.00          &                 &             &             & $\bullet$   & $\bullet$   & $\bullet$   & \multicolumn{1}{c|}{$\bullet$}           \\
\multicolumn{1}{|c|}{\textbf{28}}  & 492.32              & 98.74           & 98.68           & 98.53           & 98.32           & 98.05        	 & -               & -                  & -              & 0.22                 & -                  & 0.00          &                 &             &             &             & $\bullet$   &             & \multicolumn{1}{c|}{$\bullet$}           \\ \hline
\end{tabular}
}
\end{table}
\end{landscape}	

Looking at the Pareto sets in Table~\ref{sample1tabular} and comparing the two monotonic measures, we see that the lexicographic min-max can yield a noticeable number of additional trade-off solutions compared to the standard min-max (solutions 4, 6, 7, 8). It is likely that this difference will grow with instance size. Interestingly, some of the solutions identified with the lexicographic measure are also found by all the non-monotonic functions (solutions 4, 5, and 6), while others are not found by any of them (solutions 7, 8, and 11). In theory, the lexicographic objective also favours more equitable solutions in those cases when there exist multiple min-max optima for the same total cost. However, in our study we did not encounter such an instance.

Comparing the four relative measures, we observe that their fronts are noticeably larger than those generated with absolute measures. Although some solutions are also shared with both absolute measures (solutions 1, 2, 13, 15, 24), the Pareto sets found with relative measures are not strictly supersets of those found with absolute measures. Some solutions (e.g.~9, 12, and 23) are found only with the absolute measures, while others (e.g.~11 and 16) are found only with the relative measures. Comparing the relative measures themselves, the differences become smaller. Although in this particular instance all the solutions found by the range measure are also found by the remaining three, in general none of the relative measures is subsumed by another. We find that all the relative measures can typically be reduced to nearly 0.

With respect to the TSP-optimality of the identified Pareto-optimal solutions, one can see that TSP-optimal solutions are present at both the lower cost and lower inequality ends of all fronts, while non-TSP-optimal solutions are mainly encountered at lower inequality. Some TSP-optimal solutions (e.g.~11, 16, and 20) are identified only with non-monotonic measures, others (e.g.~9, 12, and 23) only with monotonic ones.

A closer examination of the solutions in Table~\ref{sample1tabular} reveals that inconsistency is mainly encountered at low inequality end of the fronts (solution 28 is a particularly striking example). Most of these inconsistent solutions are not TSP-optimal, and their total solution cost rapidly increases. However, even low-cost solutions may be inconsistent. Solution 3 is a typical example: one of the tours of solution 2 is lengthened while the rest remain unchanged. Comparing all the equity functions, one observes that inconsistent solutions are generated only with the non-monotonic inequality measures (in line with Theorem~\ref{T3} and Observation~\ref{T4}). Nonetheless, some TSP-optimal and consistent solutions (e.g.~solution 11) are found only with such measures. As pointed out in Observation~\ref{P1}, there exist solutions which are TSP-optimal but still inconsistent (e.g.~solutions 16 and 20).

\begin{figure}[htbp]
\FIGURE
{\includegraphics[width=0.75\textwidth]{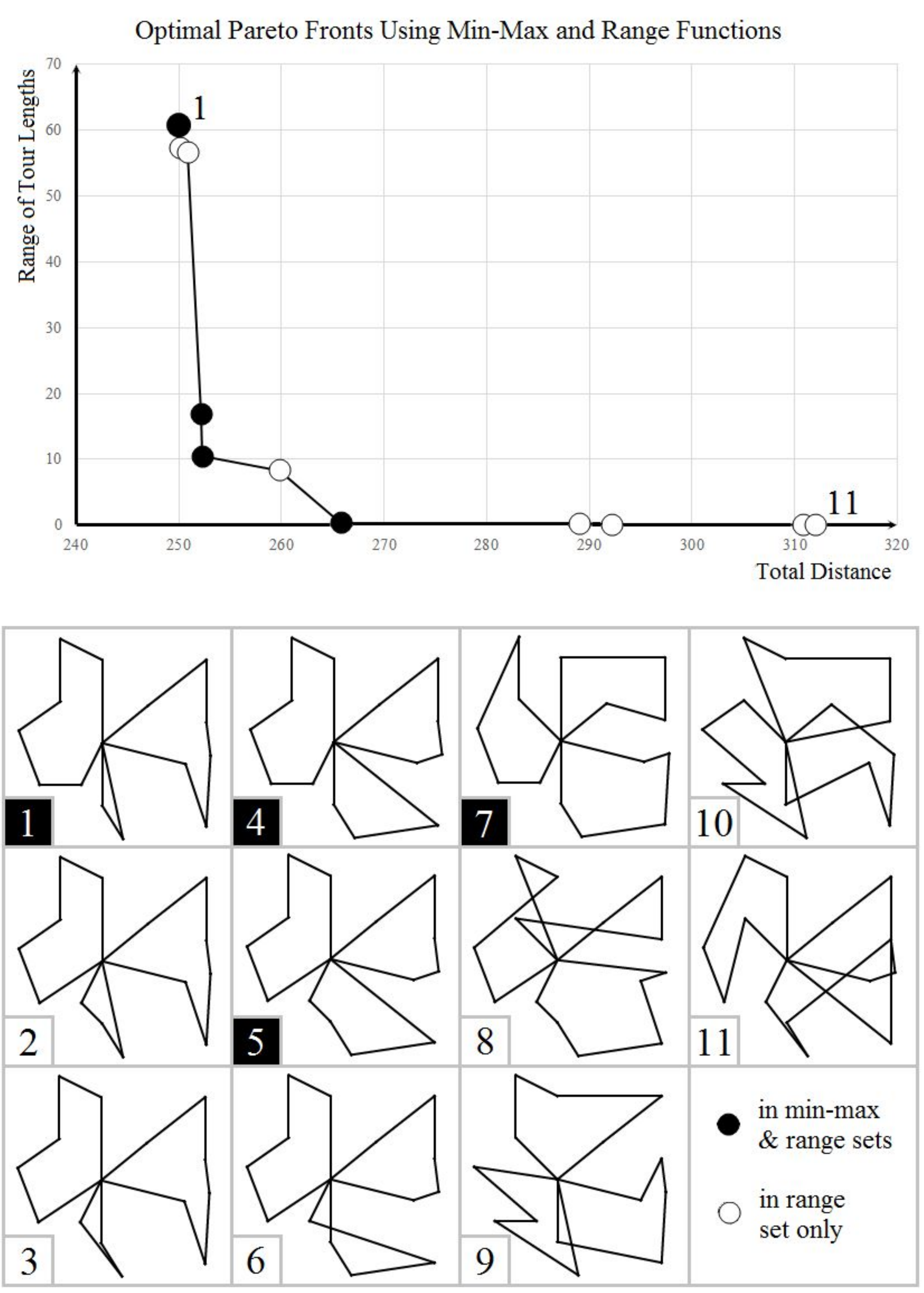}}
{Visualization of VRPRB Pareto Sets and Corresponding Solutions for a Sample Instance (\#46).\label{samplefront1}}
{}
\end{figure}

In order to show how these two aspects impact the structure of VRPRB solutions in the decision space, we visualize in Figure~\ref{samplefront1} two Pareto sets (min-max and range) and their corresponding solutions for another representative instance. In this example, the Pareto set obtained with the min-max measure is a subset of the Pareto set using the range measure. Solutions 1, 4, 5, and 9 are common to both Pareto fronts. Starting from the optimal cost solution, solution 2 is obtained by moving a customer from the medium-length tour to the shortest one. This improves the range but not the min-max, and therefore this solution is dominated in the min-max set. However, the lexicographic min-max would capture also solution 2. One can see that solution 3 is then obtained by increasing the length of the shortest tour. The resulting solution is no longer TSP-optimal, and is not consistent if solution 2 is also identified. Solutions 4 and 5 reassign the customers between the longest and shortest tours, yielding two TSP-optimal and consistent solutions. Solution 6 is a typical example of how non-monotonic measures like the range may be artificially improved by breaking the TSP-optimality of a previous solution (in this case of solution~5). This behaviour is apparent again with solutions 8, 9, 10, and especially 11, the solution with minimal range. Visual inspection reveals that, in the absence of further constraints, solutions with the latter characteristics are not likely to find acceptance in practical settings.

As mentioned earlier, the issue of non-TSP-optimal tours in the efficient sets of the VRPRB has been noted previously by \cite{jozefowiez2002}, who propose to apply a 2-opt local search to new candidate solutions prior to examining Pareto efficiency. Given the observations made above, this suggestion appears to be well-motivated. In the following section, we therefore add an even stronger TSP-optimality constraint to the model, resolve the instances, and compare how this change impacts optimal VRPRB solution sets. 

\subsection{Optimization with TSP Constraints.}

We recall that the monotonic measures (min-max and lexicographic min-max) guarantee TSP-optimality and consistency of all Pareto-efficient VRPRB solutions, as outlined in Section~\ref{implications}. Adding a TSP-optimality constraint to the model does not have any additional impact. We therefore focus our analysis in this section on the non-monotonic measures. We note that with the TSP-optimality constraint we were able to enumerate instances with up to 20 customers, but for reasons of consistency we report here the results on the same 14-customer instances used in the previous section. The results and conclusions are qualitatively the same for the larger instances.

\begin{table}[htbp]
\centering
\SingleSpacedXI
\caption{Average Number of Pareto Efficient, Previously Unidentified TSP-optimal, and Consistent Solutions Found with Non-Monotonic Equity Functions, with TSP-Optimality Constraint}
\renewcommand{\arraystretch}{1.2}
\begin{tabular}{lrcrcrcrc} \hline
                      & \multicolumn{2}{c}{\textbf{range}} & \multicolumn{2}{c}{\textbf{mean abs. dev.}} & \multicolumn{2}{c}{\textbf{standard dev.}} & \multicolumn{2}{c}{\textbf{Gini coeff.}} \\ \hline
																				& abs.	& \%	& abs.	& \%	& abs.	& \%	& abs.	& \%	\\ \hline
\textbf{Total}											& 14.87	& -			& 15.85	& -			& 16.02	& -			& 16.45	& -  		\\
\textbf{New} 												& 9.28	& 62\%	& 9.26	& 58\%	& 9.86 	& 61\%	& 10.46	& 63\%	\\
\textbf{Consistent}									& 4.50	& 30\%	& 4.52	& 28\%	& 4.75  & 29\% 	&  4.81	& 29\%	\\ \hline
\end{tabular}
\label{summary2}
\end{table}

Given the large share of inconsistent and non-TSP-optimal solutions in the previous experiment, the question arises whether and to what extent these kinds of solutions dominate solutions that would otherwise be consistent and TSP-optimal. Table~\ref{summary2} presents summary statistics for the extended model. We report again the average cardinality of optimal Pareto sets. In order to more accurately assess the impact of adding the TSP-constraint, we report the absolute and relative number of \textit{previously unidentified} (marked 'new') TSP-optimal solutions and consistent solutions. The TSP-optimal solutions identified without the TSP-constraint remain Pareto-efficient also in the extended model. To better facilitate comparisons, Table~6 lists the updated Pareto sets for the sample instance from the previous section.

\begin{landscape}
\begin{table}[htbp]
\SingleSpacedXI
\renewcommand{\arraystretch}{1.2}
\caption{Pareto Optimal Fronts per Equity Function for a Representative Instance (\#30), with TSP-Optimality Constraint} \label{sample2tabular}
\hspace*{-1cm}
\vspace*{-0.5cm}
\scalebox{0.85}{
\begin{tabular}{|cF|ccccc|FFFFFF|c|cccccc}
\hline
                                   & \textbf{}           & \multicolumn{5}{c|}{\textbf{Tour Lengths}}                                              & \multicolumn{6}{c|}{\textbf{Equity Function}}                                                                    & \textbf{}    & \multicolumn{6}{c|}{\textbf{Inconsistency}}                                                            \\ \hline
\multicolumn{1}{|c|}{\textbf{Nr.}} & \textbf{Total Cost} & \textbf{Tour 1} & \textbf{Tour 2} & \textbf{Tour 3} & \textbf{Tour 4} & \textbf{Tour 5} & \textbf{minmax (E1)} & \textbf{lex.~rank (E2)} & \textbf{range (E3)} & \textbf{MAD (E4)} & \textbf{std.~dev. (E5)} & \textbf{Gini (E6)} & \textbf{TSP} & \textbf{E1} & \textbf{E2} & \textbf{E3} & \textbf{E4} & \textbf{E5} & \multicolumn{1}{c|}{\textbf{E6}} \\ \hline
\multicolumn{1}{|c|}{\textbf{1}}	 & 340.45              & 97.15           & 78.04           & 68.96           & 53.13           & 43.17           & 97.15           & 1               & 53.98          & 15.95                & 18.92              & 0.16          & $\checkmark$    & $\circ$     & $\circ$     & $\circ$     & $\circ$     & $\circ$     & \multicolumn{1}{c|}{$\circ$}            \\
\multicolumn{1}{|c|}{\textbf{2}}   & 341.83              & 83.04           & 78.04           & 68.96           & 68.61           & 43.17           & 83.04           & 2               & 39.87          & 10.08                & 13.74              & 0.10          & $\checkmark$    & $\circ$     & $\circ$     & $\circ$     & $\circ$     & $\circ$     & \multicolumn{1}{c|}{$\circ$}            \\
\multicolumn{1}{|c|}{\textbf{3}}   & 345.43              & 83.04           & 77.81           & 72.80           & 68.61           & 43.17           & -               & 3               & -              & -                    & -                  & 0.10          & $\checkmark$    &             & $\circ$     &             &             &             & \multicolumn{1}{c|}{$\circ$}            \\
\multicolumn{1}{|c|}{\textbf{4}}   & 347.61              & 83.04           & 75.52           & 72.80           & 68.61           & 47.64           & -               & 4               & 35.40          & 9.12                 & 11.91              & 0.09          & $\checkmark$    &             & $\circ$     & $\circ$     & $\circ$     & $\circ$     & \multicolumn{1}{c|}{$\circ$}            \\
\multicolumn{1}{|c|}{\textbf{5}}   & 350.44              & 83.04           & 75.48           & 72.80           & 68.61           & 50.50           & -               & 5               & 32.54          & 8.42                 & 10.86              & 0.08          & $\checkmark$    &             & $\circ$     & $\circ$     & $\circ$     & $\circ$     & \multicolumn{1}{c|}{$\circ$}            \\
\multicolumn{1}{|c|}{\textbf{6}}   & 352.67              & 83.04           & 68.96           & 68.61           & 66.37           & 65.68           & -               & 6               & 17.36          & 5.00                 & 6.38               & 0.04          & $\checkmark$    &             & $\circ$     & $\circ$     & $\circ$     & $\circ$     & \multicolumn{1}{c|}{$\circ$}            \\
\multicolumn{1}{|c|}{\textbf{7}}   & 354.26              & 82.81           & 81.28           & 78.04           & 68.96           & 43.17           & 82.81           & 7               & -              & -                    & -                  & -             & $\checkmark$    & $\circ$     & $\circ$     &             &             &             & \multicolumn{1}{c|}{}            \\
\multicolumn{1}{|c|}{\textbf{8}}   & 355.89              & 83.04           & 69.63           & 68.96           & 68.61           & 65.65           & -               & -               & -              & 4.75                 & 6.09               & 0.04          & $\checkmark$    &             &             &             & $\circ$     & $\circ$     & \multicolumn{1}{c|}{$\circ$}            \\
\multicolumn{1}{|c|}{\textbf{9}}   & 356.07              & 81.28           & 78.04           & 74.67           & 68.96           & 53.13           & 81.28           & 8               & -              & -                    & -                  & -             & $\checkmark$    & $\circ$     & $\circ$     &             &             &             & \multicolumn{1}{c|}{}            \\
\multicolumn{1}{|c|}{\textbf{10}}  & 357.45              & 78.04           & 74.67           & 68.96           & 68.61           & 67.16           & 78.04           & 9               & 10.88          & 3.89                 & 4.16               & 0.03          & $\checkmark$    & $\circ$     & $\circ$     & $\circ$     & $\circ$     & $\circ$     & \multicolumn{1}{c|}{$\circ$}            \\
\multicolumn{1}{|c|}{\textbf{11}}  & 361.05              & 77.81           & 74.67           & 72.80           & 68.61           & 67.16           & 77.81           & 10              & 10.65          & 3.46                 & 3.90               & 0.03          & $\checkmark$    & $\circ$     & $\circ$     & $\circ$     & $\circ$     & $\circ$     & \multicolumn{1}{c|}{$\circ$}            \\
\multicolumn{1}{|c|}{\textbf{12}}  & 362.47              & 78.04           & 74.67           & 72.19           & 68.96           & 68.61           & -               & -               & 9.44           & 3.09                 & 3.56               & 0.03          & $\checkmark$    &             &             & $\bullet$   & $\bullet$   & $\bullet$   & \multicolumn{1}{c|}{$\bullet$}           \\
\multicolumn{1}{|c|}{\textbf{13}}  & 366.08              & 77.81           & 74.67           & 72.80           & 72.19           & 68.61           & -               & -               & -              & 2.42                 & 3.02               & 0.02          & $\checkmark$    &             &             &             & $\bullet$   & $\bullet$   & \multicolumn{1}{c|}{$\bullet$}           \\
\multicolumn{1}{|c|}{\textbf{14}}  & 368.96              & 77.40           & 75.48           & 74.67           & 72.80           & 68.61           & 77.40           & 11              & -              & -                    & -                  & -             & $\checkmark$    & $\circ$     & $\circ$     &             &             &             & \multicolumn{1}{c|}{}            \\
\multicolumn{1}{|c|}{\textbf{15}}  & 369.76              & 75.48           & 74.67           & 74.60           & 72.80           & 72.21           & 75.48           & 12              & 3.27           & 1.16                 & 1.24               & 0.01          & $\checkmark$    & $\circ$     & $\circ$     & $\circ$     & $\circ$     & $\circ$     & \multicolumn{1}{c|}{$\circ$}            \\
\multicolumn{1}{|c|}{\textbf{16}}  & 412.41              & 84.39           & 82.81           & 82.34           & 81.60           & 81.28           & -               & -               & 3.12           & 0.90                 & 1.10               & 0.01          & $\checkmark$    &             &             & $\bullet$   & $\bullet$   & $\bullet$   & \multicolumn{1}{c|}{$\bullet$}           \\
\multicolumn{1}{|c|}{\textbf{17}}  & 431.98              & 88.07           & 86.66           & 86.49           & 86.17           & 84.58           & -               & -               & -		          & 0.81                 & -		              & 0.01          & $\checkmark$    &             &             & 			      & $\bullet$   & 			      & \multicolumn{1}{c|}{$\bullet$}           \\
\multicolumn{1}{|c|}{\textbf{18}}  & 467.86              & 95.93           & 93.57           & 93.16           & 92.87           & 92.33           & -               & -               & -		          & -		                 & -		              & 0.01          & $\checkmark$    &             &             & 			      & 				    & 				    & \multicolumn{1}{c|}{$\bullet$}           \\
\multicolumn{1}{|c|}{\textbf{19}}  & 494.01              & 100.83          & 99.48           & 98.19           & 98.02           & 97.50           & -               & -               & -		          & -		                 & -		              & 0.01          & $\checkmark$    &             &             & 				    & 				    & 				    & \multicolumn{1}{c|}{$\bullet$}           \\
\multicolumn{1}{|c|}{\textbf{20}}  & 496.57              & 100.17          & 100.07          & 99.19           & 99.12           & 98.02           & -               & -               & 2.15           & 0.65	               & 0.78               & 0.00          & $\checkmark$    &             &             & $\bullet$   & $\bullet$   & $\bullet$   & \multicolumn{1}{c|}{$\bullet$}           \\
\multicolumn{1}{|c|}{\textbf{21}}  & 497.99              & 100.82          & 100.17          & 99.19           & 99.12           & 98.69           & -               & -               & 2.13           & -		                 & -                  & -		          & $\checkmark$    &             &             & $\bullet$   & 			      &             & \multicolumn{1}{c|}{}           \\ \hline
\end{tabular}
}
\end{table}

\vspace*{-0.2cm}

\begin{figure}[H]
\begin{center}
\includegraphics[width=.8\textwidth]{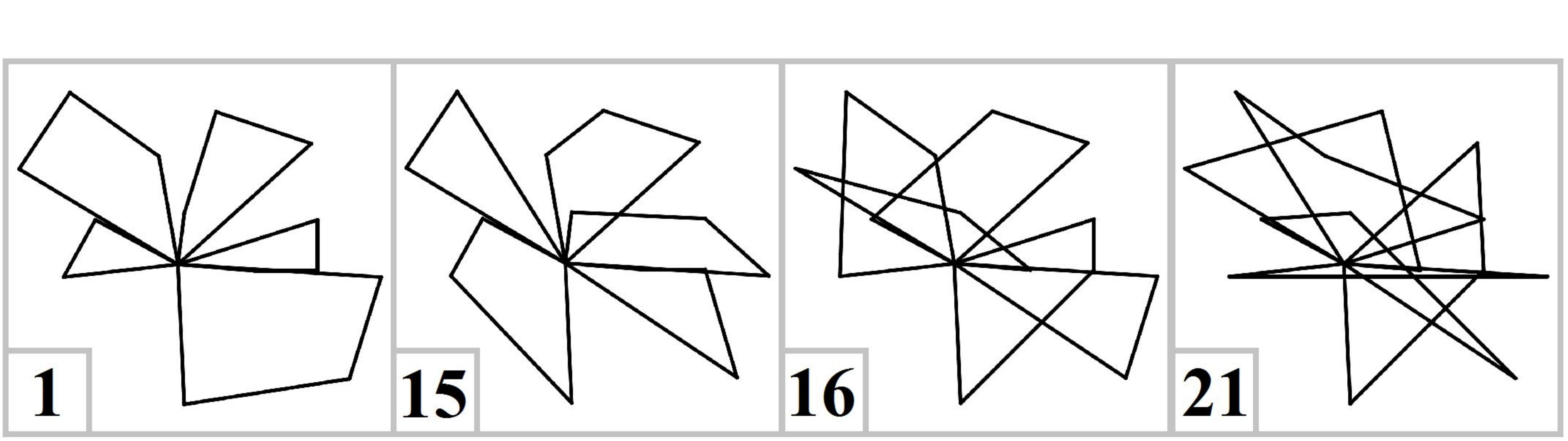}
\caption{Visualization of Selected Solutions from Table \ref{sample2tabular}} \label{solViz}
\end{center}
\end{figure}

\end{landscape}


We can see that with the TSP-constraint added, the average cardinality of fronts is approximately halved, but significantly more TSP-optimal solutions are identified in total. For all of the examined non-monotonic measures, approximately 60\% of each new Pareto set consists of previously unidentified TSP-optimal solutions. For some instances, this figure is as high as 90\%. This observation suggests that in general, a large percentage of potentially interesting solutions is dominated by non-TSP-optimal solutions in the conventional model:

\vspace*{0.2cm}
\begin{observation}
\label{P2}
Inconsistent solutions can dominate TSP-optimal solutions which would otherwise be Pareto-optimal.
\end{observation}
\vspace*{0.2cm}

However, a closer analysis of the new Pareto sets reveals that most of the solutions are still inconsistent. For all the examined non-monotonic measures, around 70\% of the solutions identified with the extended model are inconsistent. In fact, it follows from Observation~\ref{P1} that:

\vspace*{0.2cm}
\begin{observation}
\label{P3}
Tour length minimization procedures cannot eliminate inconsistent solutions from Pareto-optimal VRPRB solution sets.
\end{observation}
\vspace*{0.2cm}

Although inconsistent solutions may still appear anywhere in the Pareto set, we observe that they are most common at the high cost/low inequality area immediately after the solution that minimizes the lexicographic objective. This is not a coincidence. For the case of two vehicles (outcomes), it can be shown that all solutions with a cost higher than the lexicographic minimum are inconsistent (see Appendix~1). With a greater number of vehicles (outcomes), there is a possibility for some consistent solutions to exist beyond the lexicographic minimum. Based on our numerical study using up to five vehicles, we observe that such solutions are very rare, but do exist.

We visualize in Figure~\ref{solViz} some of the solutions from Table~\ref{sample2tabular}: the cost-optimal solution, the lexicographic minimum, the first inconsistent solution after the lexicographic minimum, and the solution with lowest inequality among all non-monotonic measures. It is difficult to imagine that a decision maker would implement solutions like solution~21, even if equity were the primary optimization objective. Numerical as well as visual inspection reveals that the equity of such solutions is artificial.

After extending the model with a TSP-constraint, all tours are TSP-optimal, but customers can still be allocated in ways that lengthen all or some TSP-optimal tours. As pointed out in Section~\ref{implications}, this is simply another way to increase the length of tours, which in turn can improve the inequality index returned by non-monotonic measures. As noted above, adding any form of tour length minimization procedure cannot solve the underlying problem caused by non-monotonicity of the equity function. With respect to the \textit{extent} of this problem, the observations made in our numerical study strongly suggest that it is not merely a special case, but rather a common feature of Pareto-optimal solutions identified by non-monotonic measures.

\section{Conclusions and Perspectives.}
\label{5-conclusion}

Classifying the existing literature as a point of departure, it is clear that issues surrounding equity in logistics are gaining increasing attention from researchers and practitioners alike. The various models and methodologies proposed by researchers show that equity considerations pose interesting theoretical challenges, and the relatively large number of application-oriented papers confirm that equity concerns are also relevant in practice. Yet despite this broad consensus about the importance and relevance of equity considerations, there has been limited discussion about just how equity should be measured and assessed in a logistics context.

\vspace*{0.2cm}
\paragraph{\bf Analytical Results.}
We began our analysis by considering six commonly applied measures of equity and comparing them. We determined that no measure satisfies every desirable axiom and none are strictly better than others in all relevant aspects. We pointed out the implications of a measure satisfying or not the Pigou-Dalton transfer principle and the monotonicity axiom. In particular, if workload is based on a variable-sum metric like tour length, then non-monotonic (relative) measures of inequality can generate Pareto sets with two undesirable consequences: Pareto-optimal solutions may contain \textbf{non-TSP-optimal} tours, and they may also be \textbf{inconsistent} with the minimization of individual workloads, i.e.~they may consist of tours whose workloads are all equal to or longer than those in other efficient solutions. Moreover, we pointed out that enforcing even TSP-optimality of tours is not sufficient for avoiding workload inconsistency.

\vspace*{0.2cm}
\paragraph{\bf Empirical Observations.}
In order to estimate the extent to which these issues may arise in practice and to assess other relevant aspects not directly implied by axiomatic properties, we performed a numerical study on small instances solvable to optimality for all examined inequality measures. We confirmed the observations made in previous studies that cost-optimal VRP solutions are usually quite poorly balanced: on average, the longest tour is found to be about twice as long as the shortest. However, the marginal cost to improve equity is usually low. In fact, with a view toward practical applications in which the cost criterion is likely to remain the primary objective, nearly 40\% of the identified efficient solutions in our study have costs less than 10\% above the optimal cost, showing that reasonable compromise solutions often exist. These observations further motivate a multi-objective approach to equity in VRPs, as well as the search for an appropriate equity measure.

Our study reveals that a conventional optimization of VRPRB instances results in a very high degree of non-TSP-optimality and inconsistency in the Pareto sets identified with non-monotonic measures. On average, less than 20\% of Pareto-optimal solutions were TSP-optimal, and less than 15\% were consistent. The literature on the VRPRB proposes to extend the model with tour length minimization procedures (e.g.~2-opt) to filter out such solutions during optimization. Unfortunately, this does not attack the underlying cause (non-monotonicity of the equity measure), and therefore cannot entirely mitigate the mentioned effects. Nonetheless, we extended our model with a TSP-constraint in order to assess the extent to which such approaches may reduce non-TSP-optimality and inconsistency in practice. We observed that around 50\% of the Pareto-optimal solutions identified with the extended model are TSP-optimal solutions which would otherwise have been dominated in the conventional approach. However, fewer than 5\% of these new solutions are consistent. We therefore conclude that tour length minimization procedures are neither sufficient in theory nor in practice to overcome the problems caused by non-monotonic inequality measures.

\vspace*{0.2cm}
\paragraph{\bf Proposed Guidelines.}
We arrive at the counter-intuitive conclusion that more sophisticated equity measures do not necessarily result in more reasonable trade-off solutions, even when combined with a min-sum cost objective. Although measures such as the range or the standard deviation are undoubtedly more informative when summarizing and comparing two workload allocations, they can lead to some paradoxes when optimized directly. This reduces the credibility of such models, which can be damaging in applications.

In light of these results, we conclude that \textbf{monotonic measures} of inequality, such as min-max or its lexicographic extension, are more appropriate for numerical optimization when workloads are based on \textbf{variable-sum metrics}. This includes the tour length metric used in the VRPRB, but also applies to the more general case of weighted aggregations of different metrics so long as at least one is variable-sum, as in some studies combining length and load balancing. The choice of the range measure for the VRPRB thus warrants re-assessment, particularly if the VRPRB is to be used as a prototypical template for equitable bi-objective VRP models.

When it comes to the choice between the two monotonic measures that we have examined, the lexicographic min-max does yield an appreciable number (some 20\% more on average) of unique trade-off solutions that are otherwise dominated using the standard min-max measure. These additional solutions tend to be in the low cost side of the Pareto sets, which is likely to be more interesting in practice. The disadvantages are that lexicographic optimization is more complicated to implement, and also to interpret since it is not immediately clear how the resulting fronts can be presented to a decision maker. Given that the size of VRP Pareto sets tends to increase exponentially, the standard min-max measure is likely a more practical choice that still identifies a reasonable absolute number of trade-off solutions.

\vspace*{0.2cm}
\paragraph{\bf Open Research Directions.}
There remain a number of promising avenues for further research on equity issues in logistics. First, there is a distinct lack of models and methods for the more general case when the assumption of anonymity/impartiality does not hold. A practically-relevant example is the distinction between full-time and part-time employees, who certainly have different levels of utility or disutility for the same absolute workload. Designing adequate models for this generalization (not necessarily restricted to logistics applications) would be a significant contribution to the equity literature, upon which the development of efficient solution methods should subsequently be based.

Second, it is worthwhile to extend existing models and methods to the multi-period case. In practice, workload considerations are commonly based on a rolling horizon perspective, especially in the case of employees' working hours. In this context, cost-minimizing solutions might appear inequitable in individual periods. Yet by allocating workloads over the horizon, it could be possible to attain an equitable overall workload distribution without significantly worsening the cost minimization objective. A relevant concern in this context would be the question to what extent balanced workloads conflict with issues such as customer and driver consistency \citep{kovacs2014} (distinct from the workload inconsistency discussed in this article).

Third, we have observed that some equitable VRP solutions lack credibility. Visual inspection reveals that some tour workloads can be improved by simple relocations or exchanges of customers, even if all tours are TSP-optimal. It is clear that trade-off solutions will usually not be optimal for any single objective, but if those objectives can be improved by simple visual inspection, then such solutions are not likely to be convincing in practice. The addition of local optimality constraints with neighbourhoods such as relocate or exchange is one possible way to ensure that trade-off solutions are credible without negatively impacting (heuristic) search methods. We are not aware of any studies examining the presence (or lack) of local optima in optimal Pareto sets of VRPs.

Finally, the theoretical literature that we have surveyed focuses almost exclusively on equity in terms of tour length, but the papers reporting on applications include also cases in which other factors are equally or more relevant. In the interest of generalizing existing theory and methodology, further research should place greater emphasis on broader definitions of workload. For example, in the small package delivery sector, the number of stops is the primary determinant of workload (see also \citealt{gulczynski2011} and \citealt{groer2009}). \cite{lin2006} and \cite{liu2006} describe applications in telecommunications and groceries delivery, in which load/demand is at least as important as tour length. For equity metrics whose sum is constant for all solutions (number of stops, demand, service time), monotonicity of the equity function is no longer a decisive property since every increase/decrease in one outcome must be accompanied by the reverse in another. This in turn implies that the PD transfer principle becomes more important, because every pair of solutions is then connected by mean-preserving transfers. It is therefore likely that for constant-sum metrics, measures such as the range will generate richer Pareto sets without the shortcomings identified in the present study. A comparative study of equity measures for this class of metrics is therefore another promising research subject.


%

\vspace*{0.3cm}
\section*{Electronic Companion}
The instances generated for this study, as well as detailed computational results for all instances, are available upon request from the authors and as part of the online version of this paper.

\vspace*{0.3cm}
\section*{Acknowledgements}
The authors extend their thanks to the two anonymous referees and the editors, whose comments and suggestions helped to improve the quality of this paper.

\vspace*{0.3cm}
\section*{Appendix}

\vspace*{0.2cm}
\begin{theorem}
\label{T5}
Let $\mathbf{x}$ be the lexicographically minimal solution to an instance of the VRPRB. For $n=2$, any solution $\mathbf{y} \neq \mathbf{x}$ with $C(\mathbf{y}) > C(\mathbf{x})$ is either inconsistent or dominated by~$\mathbf{x}$.
\end{theorem}

\proof{Proof.}
For $n=2$, $\mathbf{x} = (x_1, x_2)$, and $\mathbf{y} = (y_1, y_2)$. Let the outcomes be sorted in non-increasing order. If $\mathbf{x}$ is lexicographically minimal, then $x_1 \leq y_1$. One can distinguish the following two cases with respect to the relationship between $x_2$ and $y_2$:
\begin{itemize}[nosep]
\item If $x_1 \leq y_1$ and $x_2 \leq y_2$, then $\mathbf{y}$ is inconsistent according to Definition~\ref{def1}.
\item If $x_1 \leq y_1$ and $x_2 > y_2$, then $I(\mathbf{y}) \geq I(\mathbf{x})$ for all equity measures examined in this paper. Since $C(\mathbf{y}) > C(\mathbf{x})$ and $I(\mathbf{y}) \geq I(\mathbf{x})$, $\mathbf{y}$ is dominated by $\mathbf{x}$ and cannot be Pareto-optimal.~\Halmos
\end{itemize}
\endproof

%
%


\vspace*{0.3cm}
\bibliographystyle{informs2014trsc} 
\bibliography{bibliography} 


\end{document}